\documentclass[twocolumn,amsmath,amssymb,aps,secnumarabic,%
    nofootinbib,groupedaddress]{revtex4}
\usepackage{graphics,times,mathptmx,amsthm,pstricks}
\usepackage[small]{subfigure}
\newtheorem{theorem}{Theorem}[section]
\newtheorem{scholium}[theorem]{Scholium}
\newtheorem{corollary}[theorem]{Corollary}
\newtheorem{proposition}[theorem]{Proposition}

\theoremstyle{definition}
\theoremstyle{remark}
\newtheorem*{remark}{Remark}

\newcommand{\C}{\mathbb{C}}
\newcommand{\F}{\mathbb{F}}
\newcommand{\R}{\mathbb{R}}
\newcommand{\Z}{\mathbb{Z}}
\renewcommand{\tensor}{\otimes}

\newcommand{\ve}{\vec{e}}
\newcommand{\vp}{\vec{p}}
\newcommand{\vq}{\vec{q}}
\newcommand{\vx}{\vec{x}}
\newcommand{\vy}{\vec{y}}
\newcommand{\vz}{\vec{z}}
\newcommand{\vlam}{\vec{\lambda}}
\newcommand{\bM}{\bar{M}}
\newcommand{\bV}{\bar{V}}
\newcommand{\bz}{\bar{z}}
\newcommand{\Tr}{\mathrm{Tr}}
\newcommand{\Aut}{\mathrm{Aut}}
\newcommand{\CP}{{\C P}}
\newcommand{\RP}{{\R P}}
\newcommand{\GL}{\mathrm{GL}}
\newcommand{\SO}{\mathrm{SO}}
\newcommand{\PSU}{\mathrm{PSU}}
\newcommand{\eps}{\epsilon}
\newcommand{\defeq}{\stackrel{\mathrm{def}}{=}}
\renewcommand{\hat}{\widehat}
\renewcommand{\bar}{\overline}
\newcommand{\hT}{\hat{T}}
\newcommand{\hF}{\hat{F}}
\DeclareMathOperator{\Vol}{Vol}

\newcommand{\Spin}{\mathrm{Spin}}
\newcommand{\Ind}{\mathrm{Ind}}
\newcommand{\Inv}{\mathrm{Inv}}
\newcommand{\longto}{\longrightarrow}

\newcommand{\ie}{\textit{i.e.}}

\newcommand{\floor}[1]{\lfloor #1 \rfloor}
\newcommand{\ceil}[1]{\lceil #1 \rceil}

\newcommand{\eatline}{\vspace{-\baselineskip}}

\newcommand{\thm}[1]{Theorem~\ref{#1}}
\newcommand{\sch}[1]{Scholium~\ref{#1}}
\newcommand{\cor}[1]{Corollary~\ref{#1}}
\renewcommand{\sec}[1]{Section~\ref{#1}}

\newcommand{\prop}[1]{Proposition~\ref{#1}}

\newcommand{\fig}[1]{Figure~\ref{#1}}
\newcommand{\eq}[2]{\begin{equation}\label{#1}#2\end{equation}}

\newenvironment{fullfigure}[2]
    {\begin{figure}[htb]\begin{center}\def\fullfiga{#1}\def\fullfigb{#2}}
    {\vspace{\baselineskip}\caption{\fullfigb.}\label{\fullfiga}\end{center}\end{figure}}
\newenvironment{fullfigure*}[2]
    {\begin{figure*}[htb]\begin{center}\def\fullfiga{#1}\def\fullfigb{#2}}
    {\vspace{\baselineskip}\caption{\fullfigb.}\label{\fullfiga}\end{center}\end{figure*}}

\psset{linewidth=.5pt,dash=3pt 3pt,doublesep=.05,dotsize=1pt 5,arrowsize=2pt 3}
\newgray{gray80}{.8}
\SpecialCoor

\newcommand{\ptex}{\psline[linewidth=1pt](.1,.1)(-.1,-.1)
    \psline[linewidth=1pt](.1,-.1)(-.1,.1)}
\newcommand{\ptcir}{\pscircle[linewidth=1pt,fillstyle=solid](0,0){.12}}
\newcommand{\ptbul}{\pscircle*[linewidth=1pt](0,0){.13}}
\newcommand{\aztec}{
    \psline(0,0)(0,1)(1,1)(1,2)(2,2)(2,3)(1,3)(1,4)(0,4)(0,5)(-1,5)}
\newcommand{\afghan}{
    \psline(0,0)(1,0)(1.5,1)(2.5,1)(3,2)(2.5,3)(3,4)(2.5,5)(1.5,5)(1,6)}
\newcommand{\plu}{\psline(0,0)(1,0)(1,1)(2,1)(2,2)(1,2)(1,3)}

\begin{document}
\title{Numerical cubature from Archimedes' hat-box theorem}

\author{Greg Kuperberg}
\email{greg@math.ucdavis.edu}
\thanks{Supported by NSF grant DMS \#0306681}
\affiliation{Department of Mathematics, University of
    California, Davis, CA 95616}

\begin{abstract}
\centerline{\textit{\normalsize Dedicated to Krystyna Kuperberg on the occasion
    of her 60th birthday}}
\vspace{\baselineskip}

Archimedes' hat-box theorem states that uniform measure on a sphere projects to
uniform measure on an interval. This fact can be used to derive Simpson's
rule.  We present various constructions of, and lower bounds for, numerical
cubature formulas using moment maps as a generalization of Archimedes' theorem.
We realize some well-known cubature formulas on simplices  as projections of
spherical designs.  We combine cubature formulas on simplices and tori to make
new formulas on spheres.  In particular $S^n$ admits a $7$-cubature formula
(sometimes a $7$-design) with $O(n^4)$ points.  We establish a local lower
bound on the density of a PI cubature formula on a simplex using the moment
map.

Along the way we establish other quadrature and cubature results of independent
interest.  For each $t$, we construct a lattice trigonometric $(2t+1)$-cubature
formula in $n$ dimensions with $O(n^t)$ points.  We derive a variant of the
M\"oller lower bound using vector bundles.  And we show that Gaussian
quadrature is very sharply locally optimal among positive quadrature formulas.
\end{abstract}
\maketitle

\section{Introduction}

Let $\mu$ be a measure on $\R^n$ with finite moments. A \emph{cubature formula
of degree $t$} for $\mu$ is a set of points $F = \{\vp_a\} \subset \R^n$ and
a weight function $\vp_a \mapsto w_a \in \R$ such that
$$\int P(\vx) d\mu = P(F) \defeq \sum_{a=1}^N w_a P(\vp_a)$$
for polynomials $P$ of degree at most $t$. (If $n=1$, then $F$ is also called a
\emph{quadrature formula}.)  The formula $F$ is \emph{equal-weight} if all
$w_a$ are equal; \emph{positive} if all $w_a$ are positive; and \emph{negative}
if at least one $w_a$ is negative.  Let $X$ be the support of $\mu$.  The
formula $\F$ is \emph{interior} if every point $\vp_a$ is in the interior of
$X$; it is \emph{boundary} if every $\vp_a$ is in $X$ and some $\vp_a$ is in
$\partial X$; and otherwise it is \emph{exterior}.  We will mainly consider
positive, interior (PI) and positive, boundary (PB) cubature formulas, and we
will also assume that $\mu$ is normalized so that total measure is 1. PI
formulas are the most useful in numerical analysis \cite[Ch. 1]{Stroud:calc}.  
This application also motivates the main question of cubature formulas, which
is to determine how many points are needed for a given formula and a given
degree $t$. Equal-weight formulas that are either interior or boundary (EI or
EB) are important for other applications, in which context they are also called
\emph{$t$-designs}.

Our starting point is a connection between quadrature on the interval $[-1,1]$
and cubature on the unit sphere $S^2$, both with uniform measure.  By
Archimedes' hat-box theorem \cite{Archimedes:sphere}, the orthogonal projection
$\pi$ from $S^2$ to the $z$ coordinate preserves normalized uniform measure. In
plainer terms, for any interval $I \subset [a,b]$ or other measurable set, the
area of $\pi^{-1}(I)$ is proportional to the length of $I$; see
\fig{f:archimedes}.  (It is called the hat-box theorem because the surface area
of a hemispherical hat equals the area of the side of a cylindrical box
containing it.)  Therefore if $F$ is a $t$-cubature formula on $S^2$, its
projection $\pi(F)$ is a $t$-cubature formula on $[-1,1]$.

\begin{fullfigure}{f:archimedes}{Archimedes' hat-box theorem}
\pspicture(-2.3,-2.1)(4.85,2.1)
\psbezier[linestyle=dashed](-2,0)(-2,.425)(-1.105,.769)(0,.769)
\psbezier[linestyle=dashed](2,0)(2,.425)(1.105,.769)(0,.769)
\pscustom[fillstyle=solid,fillcolor=gray80,linestyle=none]{
    \psbezier(-1.846,0.710)(-1.846,1.102)(-1.020,1.420)(0.000,1.420)
    \psbezier(0.000,1.420)(1.020,1.420)(1.846,1.102)(1.846,0.710)
    \psarc(0,0){2}{24.624}{40.542}
    \psbezier(1.600,1.108)(1.600,1.448)(0.884,1.723)(0.000,1.723)
    \psbezier(0.000,1.723)(-0.884,1.723)(-1.600,1.448)(-1.600,1.108)
    \psarc(0,0){2}{139.458}{155.376}}
\psbezier[linestyle=dashed](-1.846,0.710)(-1.846,1.102)(-1.020,1.420)(0.000,1.420)
\psbezier[linestyle=dashed](1.846,0.710)(1.846,1.102)(1.020,1.420)(0.000,1.420)
\psbezier[linestyle=dashed](-1.600,1.108)(-1.600,1.448)(-0.884,1.723)(0.000,1.723)
\psbezier[linestyle=dashed](1.600,1.108)(1.600,1.448)(0.884,1.723)(0.000,1.723)
\pscustom[fillstyle=solid,fillcolor=gray,linestyle=none]{
    \psbezier(-1.846,0.710)(-1.846,0.318)(-1.020,0.000)(0.000,0.000)
    \psbezier(0.000,0.000)(1.020,0.000)(1.846,0.318)(1.846,0.710)
    \psarc(0,0){2}{24.624}{40.542}
    \psbezier(1.600,1.108)(1.600,0.768)(0.884,0.492)(0.000,0.492)
    \psbezier(0.000,0.492)(-0.884,0.492)(-1.600,0.768)(-1.600,1.108)
    \psarc(0,0){2}{139.458}{155.376}}
\psline(-1.846,0.710)(-1.818,0.833)\psline(1.846,0.710)(1.818,0.833)
\psline(-1.600,1.108)(-1.520,1.300)\psline(1.600,1.108)(1.520,1.300)
\psbezier(-1.846,0.710)(-1.846,0.318)(-1.020,0.000)(0.000,0.000)
\psbezier(1.846,0.710)(1.846,0.318)(1.020,0.000)(0.000,0.000)
\psbezier(-1.600,1.108)(-1.600,0.768)(-0.884,0.492)(0.000,0.492)
\psbezier(1.600,1.108)(1.600,0.768)(0.884,0.492)(0.000,0.492)
\pscircle(0,0){2}
\psbezier(-2,0)(-2,-.425)(-1.105,-.769)(0,-.769)
\psbezier(2,0)(2,-.425)(1.105,-.769)(0,-.769)
\psline{->}(2.5,0)(4,0)
\rput[b](3.25,.13){$\pi$}
\psline(4.5,-2)(4.5,2) \psline(4.4,-2)(4.6,-2) \psline(4.4,2)(4.6,2)
\psline[linewidth=2pt](4.5,.769)(4.5,1.2)
\psline(4.4,.769)(4.6,.769)\psline(4.4,1.2)(4.6,1.2)
\endpspicture
\end{fullfigure}

The 2-sphere $S^2$ has 5 especially nice cubature formulas given by the
vertices of the Platonic solids.  Their cubature properties follow purely from
a symmetry argument of Sobolev \cite{Sobolev:rotation}.  Suppose
that $G$ is the group of common symmetries of a putative cubature formula $F$
and its measure $\mu$. If $P(\vx)$ is a polynomial and $P_G(\vx)$ is the
average of its $G$-orbit, then
$$\int P_G(\vx) d\mu = \int P(\vx) d\mu \qquad P_G(F) = P(F).$$
Therefore it suffices to check $F$ for $G$-invariant polynomials. In
particular, if every $G$-invariant polynomial of degree $\le t$ is constant,
then any $G$-orbit is a $t$-design.

\begin{fullfigure*}{f:octa}{Two projections of the octahedron rule}
\subfigure[Simpson's rule]{\pspicture(-2,-2.5)(5,2.1)
\pscircle(0,0){2}
\psbezier(-2,0)(-2,-.425)(-1.105,-.769)(0,-.769)
\psbezier(2,0)(2,-.425)(1.105,-.769)(0,-.769)
\psbezier[linestyle=dashed](-2,0)(-2,.425)(-1.105,.769)(0,.769)
\psbezier[linestyle=dashed](2,0)(2,.425)(1.105,.769)(0,.769)
\rput(0.560,0.738){\ptex} \rput(-0.560,-0.738){\ptex}
\rput(1.920,-0.215){\ptex} \rput(-1.920,0.215){\ptex}
\rput(0.000,1.846){\ptex} \rput(0.000,-1.846){\ptex}
\psline{->}(2.5,0)(4,0)
\psline(4.5,-2)(4.5,2) \psline(4.4,-2)(4.6,-2) \psline(4.4,2)(4.6,2)
\rput(4.5,0){\ptex} \rput[l](4.7,0){\large $\frac23$}
\rput(4.5,2){\ptex} \rput[l](4.7,2){\large $\frac16$}
\rput(4.5,-2){\ptex} \rput[l](4.7,-2){\large $\frac16$}
\endpspicture}
\hspace{2cm}
\subfigure[2-point Gauss-Legendre rule]{\pspicture(-2,-2.5)(5,2.1)
\pscircle(0,0){2}
\psbezier(-2,0)(-2,-.425)(-1.105,-.769)(0,-.769)
\psbezier(2,0)(2,-.425)(1.105,-.769)(0,-.769)
\psbezier[linestyle=dashed](-2,0)(-2,.425)(-1.105,.769)(0,.769)
\psbezier[linestyle=dashed](2,0)(2,.425)(1.105,.769)(0,.769)
\rput(0.457,1.669){\ptex} \rput(-0.457,-1.669){\ptex}
\rput(1.129,0.612){\ptex} \rput(-1.129,-0.612){\ptex}
\rput(-1.586,0.917){\ptex} \rput(1.586,-0.917){\ptex}
\psbezier(-1.633,-1.066)(-1.633,-1.413)(-0.902,-1.694)(0.000,-1.694)
\psbezier(1.633,-1.066)(1.633,-1.413)(0.902,-1.694)(0.000,-1.694)
\psbezier[linestyle=dashed](-1.633,-1.066)(-1.633,-0.719)(-0.902,-0.438)(0.000,-0.438)
\psbezier[linestyle=dashed](1.633,-1.066)(1.633,-0.719)(0.902,-0.438)(0.000,-0.438)
\psbezier(-1.633,1.066)(-1.633,0.719)(-0.902,0.438)(0.000,0.438)
\psbezier(1.633,1.066)(1.633,0.719)(0.902,0.438)(0.000,0.438)
\psbezier[linestyle=dashed](-1.633,1.066)(-1.633,1.413)(-0.902,1.694)(0.000,1.694)
\psbezier[linestyle=dashed](1.633,1.066)(1.633,1.413)(0.902,1.694)(0.000,1.694)
\psline{->}(2.5,0)(4,0)
\psline(4.5,-2)(4.5,2) \psline(4.4,-2)(4.6,-2) \psline(4.4,2)(4.6,2)
\rput(4.5,1.1552){\ptex} \rput[l](4.7,1.155){\large $\frac12$}
\rput(4.5,-1.155){\ptex} \rput[l](4.7,-1.155){\large $\frac12$}
\endpspicture}
\eatline
\end{fullfigure*}

By Sobolev's theorem, the vertices of a regular octahedron form a 3-design on
$S^2$.  If we project this formula using Archimedes' theorem, the result is
Simpson's rule. Another projection of the same 6 points yields 2-point
Gauss-Legendre quadrature.  \fig{f:octa} shows both projections. The 8 vertices
of a cube are also a 3-design.  One projection is again 2-point Gauss-Legendre
quadrature; another is Simpson's $\frac38$ rule.  Finally the 12 vertices of a
regular icosahedron form a 5-design by symmetry.  One projection of these 12
points is 4-point Gauss-Lobatto quadrature.

The rest of this article applies toric moment maps, which generalize
Archimedes' theorem to higher dimensions, to the cubature problem.
\sec{s:projection} shows that several well-known quadrature formulas on the
interval and cubature formulas on simplices are projections of
higher-dimensional, symmetric formulas.   \sec{s:fibration} combines formulas
on tori with formulas on simplices and moment maps to make formulas on spheres
and projective spaces.  In particular, it constructs a PI 7-cubature formula on
the sphere $S^n$ with $O(n^4)$ points. Finally \sec{s:local} uses moment maps
to establish a local lower bound for the density of points in any PI cubature
formula on a simplex.  A similar lower bound holds for an arbitrary simple
convex polytope.

Along the way we establish some other quadrature and cubature results that are
not derived from moment maps but are of independent interest.   \sec{s:torus}
establishes new lattice cubature formulas on tori that are similar to cubature
formulas based on error-correcting codes \cite{Kuperberg:cubature}. In
particular it constructs, for each $t$, a trigonometric $(2t+1)$-cubature
formula on $[0,2\pi)^n$ of lattice type with $O(n^t)$ points. This improves a
construction of Cools, Novak, and Ritter with $O(n^{2t})$ points and negative
weights \cite{CNR:smolyak}, and agrees up to a $t$-dependent constant factor
with the Stroud-Mysovskikh lower bound \cite{Stroud:more,Mysovskikh:exact}. 
\sec{s:algebraic} presents a refinement of this well-known lower bound in odd
degree.  It is similar to the M\"oller bound \cite{Moller:lower}, but applies
to some new cases. \sec{s:local} also establishes that Gaussian quadrature is
very sharply locally optimal among all positive quadrature formulas
(\thm{th:sharp}).  This bound might be previously known since Gaussian
quadrature has been widely studied, but the author could not find it in the
literature.

\section{Projection constructions}
\label{s:projection}

The immediate higher-dimensional generalization of Archimedes' theorem replaces
the sphere $S^2$ by the complex manifold $\CP^n$. This manifold has a natural
metric and a natural real algebraic structure.  Concretely, assume that the
projective coordinates $(z_0:z_1:\cdots:z_n)$ of $\CP^n$ are normalized so that
$$|z_0|^2 + |z_1|^2 + \ldots + |z_n|^2 = 1.$$
Then the coordinates $z_k\bz_j$ together embed $\CP^n$ into $\C^{(n+1)^2}$ as a
real algebraic variety (with $\C^{(n+1)^2}$ interpreted as a
$2(n+1)^2$-dimensional real vector space) and a Riemannian manifold. This
embedding is familiar in quantum mechanics as the density matrix (or density
operator) formalism \cite[\S2.4]{NC:book}.  The induced metric is called the
Fubini-Study metric.  Since the metric yields a measure on $\CP^n$,
and since it is a real algebraic variety, we can consider cubature
formulas on it.

There is a projection $\pi:\CP^n \to \Delta_n$ to the $n$-simplex given by
$$\pi(z_0:z_1:\cdots:z_n) = (|z_0|^2,|z_1|^2,\ldots,|z_n|^2),$$
using normalized coordinates for $\CP^n$ and barycentric coordinates for
$\Delta_n$.  It is linear and it preserves normalized measure.  In more
abstract terms, $\pi$ has these properties because $\CP^n$ is a projective
toric variety and $\pi$ is its moment map.  Archimedes' theorem is a
description of the moment map of $\CP^1 \cong S^2$.  Thus, if $F$ is an
interior $t$-cubature formula on $\CP^n$, then $\pi(F)$ is a $t$-cubature
formula on $\Delta_n$.

\begin{fullfigure}{f:trisimp}
    {A 2-dimensional generalization of Simpson's rule}
\pspicture(-2,-1.25)(2,2.25)
\pspolygon(0,2)(-1.732,-1)(1.732,-1)
\rput(0,2){\ptex}
\rput(-1.732,-1){\ptex}
\rput(1.732,-1){\ptex}
\rput(0,0){\ptcir}
\rput[l](.2,2){\large $\frac1{12}$}
\rput[l](.2,0){\large $\frac34$}
\endpspicture\eatline
\end{fullfigure}

Ivanovi\'c, Wootters, and Fields \cite{Ivanovic:formal,WF:unbiased} defined one
interesting family of 2-designs on $\CP^{q-1}$ for $q = p^k$ a prime power.  If
$p$ is odd, then the 2-design is the orbit of a standard basis vector $e_k$ in
the group generated by cyclic permutation and linear operators of the form
$$L(e_k) = \omega^{\Tr_p(ak^2+bk+c)}e_k,$$
where $\omega$ is a $p$th root of unity and $\Tr_p$ is the $\F_p$ trace
function on $\F_q$.  The construction is more complicated when $p=2$.  In
either case, the standard basis projects to the vertices of $\Delta_{q-1}$ and
the other $q^2$ vectors project to the center.  The result is a standard degree
2 generalization of Simpson's rule for $\Delta_{q-1}$, shown in \fig{f:trisimp}
when $q=3$.

Other interesting designs and cubature formulas on $\CP^{n-1}$ come from
designs and formulas on $S^{2n-1}$. The generalized Hopf fibration
$$h:S^{2n-1} \to \CP^{n-1}$$
is a quadratic, volume-preserving map from $S^{2n-1}$ to $\CP^{n-1}$. Namely,
if we place $S^{2n-1}$ in $\C^n$, $h$ takes each point to the complex line
containing it. The map $h$ projects a $2t$- or $2t+1$-cubature formula on
$S^{2n-1}$ to a $t$-cubature formula on $\CP^{n-1}$.

\begin{fullfigure*}{f:tetra}
    {3-cubature formulas on $\Delta_3$ from the $E_8$ root system}
\subfigure[8 points]{\pspicture(-3,-2)(3,2.5)
\psline[linestyle=dashed](0.000,2.769)(0.621,0.138)
\psline[linestyle=dashed](2.079,-1.661)(0.621,0.138)
\psline[linestyle=dashed](-2.700,-1.247)(0.621,0.138)
\psline(0.000,2.769)(2.079,-1.661)
\psline(0.000,2.769)(-2.700,-1.247)
\psline(2.079,-1.661)(-2.700,-1.247)
\rput(0.000,2.769){\ptex} \rput(0.621,0.138){\ptex} \rput(2.079,-1.661){\ptex}
\rput(-2.700,-1.247){\ptex} \rput(0.000,-0.923){\ptcir}
\rput(-0.103,0.131){\ptcir} \rput(-0.693,0.554){\ptcir} \rput(0.900,0.416){\ptcir}
\rput[l](0.2,2.769){\large $\frac1{40}$}
\rput[l](0.2,-0.923){\large $\frac9{40}$}
\endpspicture}
\hspace{2cm}
\subfigure[11 points]{\pspicture(-3,-2)(3,2.5)
\psline[linestyle=dashed](0.000,2.769)(0.621,0.138)
\psline[linestyle=dashed](2.079,-1.661)(0.621,0.138)
\psline[linestyle=dashed](-2.700,-1.247)(0.621,0.138)
\psline(0.000,2.769)(2.079,-1.661)
\psline(0.000,2.769)(-2.700,-1.247)
\psline(2.079,-1.661)(-2.700,-1.247)
\rput(0.000,2.769){\ptex} \rput(0.621,0.138){\ptex} \rput(2.079,-1.661){\ptex}
\rput(-2.700,-1.247){\ptex} \rput(0.310,1.454){\ptcir} \rput(1.040,0.554){\ptcir}
\rput(-1.350,0.761){\ptcir} \rput(1.350,-0.761){\ptcir}
\rput(-1.040,-0.554){\ptcir} \rput(-0.310,-1.454){\ptcir}
\rput(0.000,0.000){\ptbul}
\rput[l](0.2,2.769){\large $\frac1{60}$}
\rput[br](-.15,.1){\large $\frac{32}{60}$}
\rput[l](1.240,0.554){\large $\frac{4}{60}$}
\endpspicture}\eatline
\end{fullfigure*}

One interesting example is the 240 roots of the $E_8$ root system, which are a
$7$-design as well as the solution to the sphere kissing problem in $\R^8$
\cite[\S14.2]{CS:splag}.  The root system has two natural positions in
$\C^4$.   In the first position, it is generated from the two points
$$(1,1,1,0) \qquad (1-\omega,0,0,0)$$
by freely permuting the first 3 coordinates, applying the map $(a,b,c,d)
\mapsto (d,a,-b,c)$, and multiplying any one coordinate by $\omega$, a cube
root of unity.   In the second position, it is generated from the three points
$$(1,1,1,1) \qquad (2,0,0,0) \qquad (1+i,1+i,0,0)$$
by freely permuting the four coordinates and multiplying any two coordinates by
$i$.   These two positions respectively exhibit the Eisenstein and Gaussian
lattice structures of the $E_8$ lattice.  The Hopf fibration sends the
Eisenstein position of the root system to a 40-point 3-design in $\CP^3$ and
the Gaussian position to a 60-point 3-design.  Then the moment map projects
these two 3-designs to 3-cubature formulas for the tetrahedron $\Delta_3$ that
appear in Abramowitz and Stegun \cite[p. 895]{AS:handbook}. They have 8 and 11
points, respectively, and are shown in \fig{f:tetra}.

The composition $\pi \circ h$ of the moment map and the Hopf fibration is a
torus fibration $\tau_2:\R^{2n} \to \Delta_{n-1}$ that does not fully depend on
the complex structure $\R^{2n} = \C^n$, but only on the decomposition of
$\R^{2n}$ into $n$ orthogonal planes.  Explicitly, the map is:
$$\tau_2(x_1,\ldots,x_{2n}) =
    (x_1^2 + x_2^2,x_3^2 + x_4^2,\ldots,x_{2n-1}^2 + x_{2n}^2).$$
This projection is analogous to a map $\tau_1:S^{n-1} \to \Delta_{n-1}$ defined
by Xu \cite{Xu:simplices}:
$$\tau_1(x_1,\ldots,x_n) = (x_1^2,\ldots,x_n^2).$$
The Xu map does not preserve uniform measure.  Rather, it takes
uniform measure on the sphere to the measure with weight function
$$w_1(\vy) = \frac{2^n\pi^{n/2}}{\frac{n}{2}!n\sqrt{y_0y_1y_2\ldots y_{n-1}}}$$
in barycentric coordinates.

In the case of the $E_8$ root system, one interesting set of orthogonal planes
is the 4 eigenplanes of the abelian subgroup of $\Aut(E_8)$ of the form $C_5
\times C_5$.  Eric Rains \cite{Rains:personal} has computed the corresponding
3-cubature formula on $\Delta_3$ using Magma \cite{Magma}.  In barycentric
coordinates on $\Delta_3$, its points and weights are the orbits of the two
weighted points
\begin{align*}
\vp_1 &= \frac1{10}(0,0,5-\sqrt{5},5+\sqrt{5}) & w_1 = \frac1{24} \\
\vp_2 &= \frac1{10}(2,2,3+\sqrt{5},3-\sqrt{5}) & w_2 = \frac5{24}
\end{align*}
under the action of the coordinate permutations $(3 4)$ and $(1 3)(2 4)$.  In
particular, it has 8 points. In conclusion, at least three interesting
3-cubature formulas for $\Delta_3$ arise as projections of $E_8$ root system. 
The root system model explains the simple rational values of the weights.

The $E_8$ lattice is one of four widely studied and highly symmetric lattices in
low dimensions; the other three are the Coxeter-Todd lattice $K_{12}$ in
$\R^{12}$, the Barnes-Wall lattice $\Lambda_{16}$ in $\R^{16}$, and the Leech
lattice $\Lambda_{24}$ in $\R^{24}$ \cite[Ch. 4]{CS:splag}.  In each case,
cases, the set of short vectors has transitive symmetry, and in each case,
Sobolev's theorem establishes its degree as a spherical design.

The 756 short vectors of $K_{12}$ form a $7$-design on $S^{11}$.  In one of its
several presentations as an Eisenstein lattice in $\C^6$ (the ``3-base''
presentation \cite[\S7.8]{CS:splag}), the short vectors are generated from the
two points
$$(1,1,1,1,1,1)\qquad (1-\omega,\omega-1,0,0,0,0)$$
by freely permuting coordinates, multiplying the coordinates by powers of
$\omega$ whose exponents sum to $0$, and negating all coordinates.  The
projection $\tau_2$ sends these points to a 16-point $3$-cubature formula on
$\Delta_5$ generated from the points
\begin{align*}
\vp_1 &= \frac12(1,1,0,0,0,0) & w_1 &= \frac1{42} \\
\vp_2 &= \frac16(1,1,1,1,1,1) & w_2 &= \frac{27}{42}
\end{align*}
by freely permuting coordinates. This formula was found by 
Stroud \cite{Stroud:some3,Stroud:calc}.

The 4320 short vectors of $\Lambda_{16}$ form a $7$-design on $S^{15}$. In its
simplest position (which exhibits its Gaussian lattice structure), the short
vectors are generated from the two vectors
\begin{gather*}
(1,1,1,1,1,1,1,1,0,0,0,0,0,0,0,0) \\ (2,2,0,0,0,0,0,0,0,0,0,0,0,0,0,0)
\end{gather*}
by permuting coordinates under the group $\GL(4,2) \ltimes (\Z/2)^4$ of affine
automorphisms of $(\Z/2)^4$, together with sign changes that keep the
coordinate sums divisible by 4.  The projection $\tau_2$ sends these points
to a 51-point $3$-cubature formula on $\Delta_7$ generated from the points
\begin{align*}
\vp_1 &=        (1,0,0,0,0,0,0,0) & w_1 &= \frac1{1080} \\
\vp_2 &= \frac12(1,1,0,0,0,0,0,0) & w_2 &= \frac1{270} \\
\vp_3 &= \frac14(1,1,1,1,0,0,0,0) & w_3 &= \frac4{135} \\
\vp_4 &= \frac18(1,1,1,1,1,1,1,1) & w_4 &= \frac{64}{135}.
\end{align*}
under the action of the affine group $\GL(3,2) \ltimes (\Z/2)^3$.  This is not
an optimal PI 3-cubature formula, because the orbit of $\vp_2$ can be
eliminated, leaving only 23 points. But it does have a novel property:  Instead
of full symmetrization, the orbit of $\vp_3$ is in the pattern of the (8,4,3)
Steiner system. But this is as good as full symmetrization for 3-cubature,
because any monomial of degree 3 involves at most 3 coordinates. The structure
of this Barnes-Wall projection led the author to relate cubature to
combinatorial $t$-designs and orthogonal arrays \cite{Kuperberg:cubature}.

The above position of $\Lambda_{16}$ is compatible with its Gaussian lattice
structure.   Eric Rains found another interesting position which is compatible
with an Eisenstein lattice structure.  The corresponding 3-cubature formula on
$\Delta_7$ has 50 points.  They are generated from
\begin{align*}
\vp_1 &=           (1,0,0,0,0,0,0,0) & w_1 &= \frac1{720} \\
\vp_2 &=    \frac14(1,1,1,1,0,0,0,0) & w_2 &= \frac1{90} \\
\vp_3 &=    \frac13(1,1,0,0,1,0,0,0) & w_3 &= \frac1{80} \\
\vp_4 &= \frac1{12}(4,0,4,0,1,1,3,3) & w_4 &= \frac1{60} \\
\vp_5 &= \frac1{12}(4,0,4,0,1,1,1,1) & w_5 &= \frac1{40} \\
\vp_6 &= \frac1{12}(3,1,3,1,1,1,1,1) & w_6 &= \frac1{30} \\
\vp_7 &= \frac1{12}(3,1,1,1,1,1,3,1) & w_7 &= \frac1{30}
\end{align*}
by the coordinate permutations $(1 2)$, $(1 3)(2 4)(5 7)(6 8)$,
and $(1 5)(2 6)(3 7)(4 8)$.

The 196560 short vectors of the Leech lattice form a $11$-design on $S^{23}$.  
The lattice has a space Eisenstein lattice structure which Conway and Sloane
call the \emph{complex Leech lattice} \cite[\S7.8]{CS:splag}.  The complex
basis that they give leads to a 5-cubature formula on $\Delta_{11}$ generated
by the points
\begin{align*}
\vp_1 &=    \frac12(1,1,0,0,0,0,0,0,0,0,0,0) & w_1 &= \frac1{10920} \\
\vp_2 &=    \frac16(1,1,1,1,1,1,0,0,0,0,0,0) & w_2 &= \frac9{3640} \\
\vp_3 &= \frac1{18}(7,1,1,1,1,1,1,1,1,1,1,1) & w_3 &= \frac{27}{1820} \\
\vp_4 &= \frac1{18}(4,4,1,1,1,1,1,1,1,1,1,1) & w_4 &= \frac{27}{3640}
\end{align*}
by the action of the Mathieu group $M_{12}$.  In other words, the coordinates
of $\vp_2$ are permuted in the pattern of the (12,6,5) Steiner system and the
points of the other coordinate are permuted freely.  The total is 276 points.
Another interesting basis of plane consists of the mutual eigenplanes of the
$(\Z/5)^3$ subgroup of the isometry group of the Leech lattice.  Eric Rains has
computed that the corresponding 5-cubature formula on $\Delta_{11}$ has 498
points, consisting of 22 orbits of the surviving coordinate permutations.
However, since none of the Barnes-Wall formulas on $\Delta_7$ are optimal, it
is not clear that the smaller of these formulas is either.

\section{Torus constructions}
\label{s:torus}

The constructions in the next section depend on an auxiliary case that
generally works out better than cubature on $S^n$, $\CP^n$, and $\Delta_n$,
namely cubature on algebraic tori.  There is a developed theory for a special
case of this problem known as trigonometric cubature
\cite{CS:minimal,CL:moderate}.  We will describe a more general class of
problems, with one new result for the classic trigonometric cubature problem
(\thm{th:craig}).

Consider a torus group $T \cong (S^1)^n$ together with a faithful linear action
on some real vector space $V \cong \R^N$.  Then we can identify $T$ with any
faithful orbit $O$ to give it a real algebraic structure.   Since $T$ is a
compact group, it also comes with Haar measure (\ie, uniform measure).  Given
both structures, we can then consider cubature formulas for $T$. If a cubature
formula $F$ is a $t$-design and forms a subgroup of $T$, then it is called a
\emph{lattice formula}, or an \emph{additive $t$-design}.

\begin{proposition} The lattice cubature problem on $T$ is
equivalent to a lattice packing problem as follows:
\begin{description}
\item[1.] The real algebraic structure on $T$ does not depend on the orbit $O$
or the base point chosen on $O$.  The ring of polynomials on $T$ is the same as
the character ring $R(T)$.
\item[2.] Every character $\chi:T \to \C$ is homogeneous as a polynomial on
$T$.  Its degree defines a norm on $\hT$, the character group of $T$.  The norm
is generated by unit steps corresponding to the characters that appear in $V
\tensor \C$.
\item[3.] The characters that are constant on a subgroup $F \subset T$ form a
sublattice $\hF \subseteq \hT$.  This correspondence is
a bijection between finite subgroups
and sublattices such that $|F| = [\hT:\hF]$.
\item[4.] The subgroup $F$ is a $t$-design if and only if $\hF$ has
minimum distance $d = t+1$.
\end{description}
\label{p:torus}
\end{proposition}

\begin{fullfigure*}{f:poly}{Polyomino and polyhex tilings that lead to lattice rules}
\subfigure[Aztec diamonds for $T(\SO(4))$]{\psset{unit=.4cm}\pspicture*(.5,.5)(15.5,12.5)
\psframe(.5,.5)(15.5,12.5)
\multiput(0,0)(-1,5){2}{\aztec}\multiput(3,-2)(-1,5){3}{\aztec}
\multiput(6,-4)(-1,5){4}{\aztec}\multiput(6,-4)(-1,5){4}{\aztec}
\multiput(8,-1)(-1,5){3}{\aztec}\multiput(11,-3)(-1,5){4}{\aztec}
\multiput(13,0)(-1,5){3}{\aztec}\multiput(16,-2)(-1,5){3}{\aztec}
\endpspicture}
\hspace{1cm}
\subfigure[Afghan hexagons for $T(\PSU(3))$]{\pspicture*(.2,.2)(6.2,5)\psframe(.2,.2)(6.2,5)
\psset{xunit=.4cm,yunit=.3464cm}
\rput(1.5,-5){\afghan}\rput(0,0){\afghan}\rput(-1.5,5){\afghan}\rput(-3,10){\afghan}
\rput(4.5,-1){\afghan}\rput(3,4){\afghan}\rput(1.5,9){\afghan}\rput(0,14){\afghan}
\rput(9,-2){\afghan}\rput(7.5,3){\afghan}\rput(6,8){\afghan}\rput(4.5,13){\afghan}
\rput(13.5,-3){\afghan}\rput(12,2){\afghan}\rput(10.5,7){\afghan}\rput(9,12){\afghan}
\rput(15,6){\afghan}\rput(13.5,11){\afghan}
\endpspicture}
\eatline
\end{fullfigure*}

The proof of \prop{p:torus} is lengthy but routine and can be left as an
exercise for the reader.  It is essentially established in the literature when
$T = T(\SO(2n))$ acts on $\R^{2n}$ by separate rotations in $n$ orthogonal
planes.  This case is equivalent to the (cubic) trigonometric cubature problem,
defined as cubature formulas on the $n$-cube $[0,2\pi)^n$ which are exact for
trigonometric polynomials of degree $t$ \cite{CS:minimal}.  All of the
arguments generalize without change.

When $T = T(\SO(2n))$, $\hT$ is naturally identified with $\Z^n$, and its norm
is the $\ell_1$ or taxicab norm.  Another torus of interest to us  is $T =
T(\PSU(n+1))$, the group of diagonal unitary matrices with determinant 1 modulo
its center.  It acts on $\C^{(n+1)^2}$, interpreted as the space of $(n+1)
\times (n+1)$ complex matrices, by conjugation.  In this case $\hT = A_n$, the
root lattice  of $\PSU(n+1)$, and its norm is defined by taking the roots of
$A_n$ as unit steps.

\begin{theorem} Given a real algebraic torus $T$ of dimension $n$, let $K
\subset \hT \tensor \R$ be the real convex hull of the unit steps in $\hT$.
Let $\delta_L(K)$ be the lattice packing density of $K$, and let
$\Vol K$ be the volume of $K$ normalized by $\hT$. Let $t \ge 0$
and let $d = t+1$.
Then the best additive
$t$-design $F$ on $T$ has at least
$$\frac{d^n(\Vol K)}{2^n\delta_L(K)} \le |F| \le
\frac{d^n(\Vol K)}{2^n\delta_L(K)}(1 + O(t^{-1}))$$
points.
\label{th:hight}
\end{theorem}

\thm{th:hight} has been noted independently by several people for trigonometric
cubature, but may originally be due to Frolov \cite{Frolov:lattice}.  In
outline, a lattice $\hF \subset \hT$ with minimum distance $d$ produces a
packing of the dilated body $\frac{d}2K$.  The packing density $\delta_L(K)$
then yields a lower bound on the index of $\hF$. On the other hand, if
$\Lambda$ is the center lattice of the best packing of $K$, then when $t$ is
large, $\frac{d}2\Lambda$ can be approximated by a sublattice of $\hT$.  This
establishes the upper bound.

Note also that the best $\Lambda$ has rational coordinates relative to $\hT$
(or they can be made rational if $\Lambda$ is not unique), because $K$ is a
rational polytope.  Thus there exist special distances $d$ such that the best
$F$ has exactly
$$\frac{d^n(\Vol K)}{2^n\delta_L(K)}$$
points.  Also if some $d$ achieves exactitude, then so does $kd$ for every $k >
1$.

If $T = T(\SO(2n))$ is the standard cubic $n$-torus, then $K$ is the $n$-cross
polytope $C_n^*$.   For example, Minkowski established that the lattice packing
density of the regular octahedron $C_3^*$ is $\frac{18}{19}$.  So there
exists an additive $5$-design on $T(\SO(6))$ with 38 points
\cite{Frolov:lattice,Noskov:periodic}.

Since $C_2^*$ is a square, its packing density is 1. Noskov
\cite{Noskov:periodic} found the best discrete approximation to this packing
for every distance $d$ to obtain lattice rules for $T(\SO(4))$.  If $d = 2s$,
then the best approximation is exact and there is a $(2s-1)$-design with $2s^2$
points.  If $d = 2s+1$, then the best approximation  corresponds to the tiling
of $\Z^2$ by the discrete $\ell_1$ ball of radius $s$, or the tiling of the
plane $\R^2$ by certain Aztec diamonds, as shown in \fig{f:poly}(a). The ball
and the corresponding $2s$-design have $s^2+(s+1)^2$ points.

Noskov's designs have a counterpart for $T(\PSU(3))$,  where $\hat{T(\PSU(3))}
= A_2$ is the triangular lattice. If we identify $A_2$ with the Eisenstein
integers $\Z[\omega]$, then the highest-density lattice with minimum distance
$d$ is the ideal generated by
$$\lfloor \frac{d}{2} \rfloor - \omega\lfloor \frac{d+1}{2} \rfloor.$$
When $d=2s$, the dual $(d-1)$-design has $3s^2$ points and
exactly matches the tiling of the plane by regular hexagons.
When $d=2s+1$, it has $3s^2+3s+1$ points and corresponds
to a tiling the plane by the hexagonal polyhex of order
$s$ (an ``Afghan hexagon''), as shown in \fig{f:poly}(b).

\begin{theorem} Let $t \ge 0$. The torus $T(\SO(2n))$ has a $(2t+1)$-design
with $O(n^t)$ points.  More precisely it has a $2t$-design with
$(2n)^t(1+o(1))$ points as $n \to \infty$ and a $(2t+1)$-design with twice as
many points. The torus $T(\PSU(n+1))$ has a $t$-design with $n^t(1+o(1))$
points as $n \to \infty$.
\label{th:craig}
\end{theorem}

\begin{remark} \thm{th:craig} can be compared with a prior result by Cools, Novak, and Ritter
\cite{CNR:smolyak}, who obtained NI $t$-cubature formulas with $O(n^t)$ points
and negative weights.  Another comparison is with the lower bound due to Stroud
and Mysovskikh \cite{Stroud:more,Mysovskikh:exact} for trigonometric
$2t$-cubature:
$$|F| \ge \frac{(2n)^t(1-o(1))}{t!}.$$
The M\"oller bound applies to trigonometric $(2t+1)$-cubature
in its interpretation as cubature on $T(\SO(2n))$ because it is
a centrally symmetric algebraic variety.  It yields:
$$|F| \ge \frac{2(2n)^t(1-o(1))}{t!}.$$
\sec{s:algebraic} establishes an analogous lower bound
for $t$-cubature on $T(\PSU(n))$ (\cor{c:molpsu}):
$$|F| \ge \frac{n^t(1-o(1))}{\ceil{t/2}!\floor{t/2}!}.$$

Thus for each $t$, \thm{th:craig} is asymptotically optimal to within a
constant factor, even though the lower bounds do not require $F$ to be
positive or interior.
\end{remark}

\begin{proof} By \prop{p:torus}, our task is to find suitable lattices in $\Z^n
= \hat{T(\SO(2n))}$ and $A_n = \hat{T(\PSU(n+1))}$.  Our task is fulfilled by
Craig lattices \cite[\S8.6]{CS:splag} in $A_n$ and skew analogues of Craig
lattices in $\Z^n$.  We describe the $A_n$ case first.

We can model $A_n$ as the set of points in $\Z^{n+1}$ with zero coordinate sum.
Let $p \ge n+1$ be prime, and index the standard basis $\{\ve_a\}$ of
$\Z^{n+1}$ by some subset $N \subset \Z/p$.  Let $p > t > 0$, and define a
linear map $\phi:\Z^{n+1} \to (\Z/p)^t$ by
$$\phi(\ve_a) = (a,a^2,\ldots,a^t).$$
We define the lattice
$$\Lambda^{(t)}(A_n) = \ker \phi \cap A_n.$$
Plainly the index of $\Lambda^{(t)}(A_n)$ is at most $p^t = n^t(1+o(1))$. (If
$n$ is large and $p \approx n$, it is $p^t$, because any lower power of $p$
would violate the Stroud-Mysovskikh bound.)

We claim that the distance of $\Lambda^{(t)}(A_n)$ is $t+1$.  To show this, we
will show that $\phi$ is injective on the simplex $\Delta_n^{(t)} \subset
\Z_{\ge 0}^{n+1}$ of non-negative vectors with coordinate sum $t$.  A vector
$\vx \in A_n$ with root-step length at most $t$ can be expressed as the
difference of two vectors in $\Delta_n^{(t)}$; therefore injectivity shows that
none of these vectors lie in $\Lambda^{(t)}(A_n)$.

We can interpret a vector $\vx \in \Delta_n^{(t)}$
as a multiset of $S$ over the set $\{0,\ldots,n\}$ with $|S| = t$: if
$$\vx = \sum_a m_a \ve_a,$$
then $m_a$ is the multiplicity of $a \in S$.  In this interpretation,
$\phi(\vx)$ is the list of power sums
$$\sum_{a \in S} a^k$$
for $1 \le k \le t$.  By standard inversion formulas \cite{Stanley:enumerative2},
these power sums determine the elementary symmetric functions of the
elements of $S$ when $p > t$, which are the coefficients of the
polynomial
$$\prod_{a \in S} (x-a).$$
Thus $\phi(\vx)$ determines $S$ as a multiset and the vector $\vx$,
and it is injective on $\Delta_n^{(t)}$.

For $\Z^n$ (with the $\ell_1$ norm), let $p > 2n$ be prime. Index the
standard basis $\{\ve_a\}$  of $\Z^n$ by some subset $N \subset \Z/p$ such that
$N$ is disjoint from $-N$. Define $\phi:\Z^n \to (\Z/p)^t$ by
$$\hat{\phi}(\ve_a) = (a,a^3,a^5,\ldots,a^{2t-1}),$$
and define
$$\Lambda^{(t)}(\Z^n) = \ker \hat{\phi}.$$
Then the index of $\Lambda^{(t)}(\Z_n)$ is again at most (and usually exactly)
$p^t = n^1(1+o(1))$. Its distance property can be explained by embedding $\Z^n$
isometrically into $A_{2n}$ using the map
$$\alpha:\ve_a \mapsto \ve_a - \ve_{-a}.$$
Then
$$\Lambda^{(t)}(\Z^n) = \alpha^{-1}(\Lambda^{(2t)}(A_n)).$$
Since $\Lambda^{(2t)}(A_n)$ has distance at least $2t+1$, so does
$\Lambda^{(t)}(\Z^n)$. We can
boost the distance to $2t+2$ by passing to its even-sum sublattice.
\end{proof}

\begin{fullfigure}{f:plus}{A Hamming-like ``plus'' tiling}
\psset{unit=.4cm}\pspicture*(.5,.5)(15.5,10.5)
\psframe(.5,.5)(15.5,10.5)
\multiput(1,-2)(-1,2){3}{\plu}\multiput(3,-1)(-1,2){5}{\plu}
\multiput(6,-2)(-1,2){7}{\plu}\multiput(8,-1)(-1,2){6}{\plu}
\multiput(11,-2)(-1,2){7}{\plu}\multiput(13,-1)(-1,2){6}{\plu}
\multiput(14,2)(-1,2){5}{\plu}\multiput(15,5)(-1,2){3}{\plu}
\rput(15,10){\plu}
\endpspicture
\end{fullfigure}

\begin{remark} When $t=1$, the number $p$ in the proof of \thm{th:craig} need
not be prime, and the lattices $\Lambda^{(1)}(\Z_n)$ and $\Lambda^{(1)}(A_n)$
produce lattice tilings of the ball of $\ell_1$-radius 1 in $\Z^n$ and the
combinatorial simplex $\Lambda^{(1)}$ in $A_n$. For example, when $n=2$, they
are equivalent to familiar tilings of the plus pentomino (\fig{f:plus}) and the
triangle trihex. The plus tiling resembles combinatorial tilings coming from
Hamming codes \cite[\S3.2]{CS:splag}. More generally, Craig lattices resemble
low-distance BCH codes.  This resemblance is what led the author to
\thm{th:craig}.
\end{remark}

\section{Fibration constructions}
\label{s:fibration}

The projection construction in \sec{s:projection}, while instructive, is
backwards in a sense:  It is harder to make $t$-cubature formulas for
$\CP^{n-1}$ and $S^{2n-1}$ than for $\Delta_{n-1}$ for most values of $n$ and
$t$.  In this section we will use the same projections to lift cubature
formulas to spheres and projective spaces from simplices. The construction also
requires the definition and constructions of cubature formulas on tori from
\sec{s:torus}.

\begin{theorem} Let $\alpha:X \to Y$ be one of the three projections $h$, $\pi$, or
$\tau_2$, and let $T$ be a generic fiber.  Let $s = 2t+1$ when $X = S^{2n-1}$
and $s=t$ when $X = \CP^{n-1}$. Given an interior (or boundary) $t$-cubature
formula $F$ for $Y$ and an interior $s$-cubature formula $F_T$ for $T$,
there is a twisted product $s$-cubature formula $F_X = F_T \ltimes F_Y$
for $X$.  It satisfies $|F_X| = |F_T|\;|F_Y|$
and it inherits positivity from its factors.  In the boundary
case, $|F_X| \le |F_T|\;|F_Y|$.
\label{th:fiber} \end{theorem}

Note that in the three cases, $T$ is isomorphic to $S^1$,
$T(\SO(2n))$, and $T(\PSU(n))$, respectively.

\begin{proof}
Let $\sigma_Y$ be the discrete measure on $Y$ corresponding to the cubature
formula $F_Y$, and let $\sigma_X = \alpha^*(\sigma_Y)$ be the pull-back of
$\sigma_Y$ to $X$.  In other words, for each point $p$ of weight $w$ in
$F_Y$, $\sigma_X$ has a term consisting of uniform measure on the torus fiber
$\alpha^{-1}(p)$.  Also let $\mu_X$ and $\mu_Y$ be normalized uniform measure
on $X$ and $Y$.

We claim that
$$\int_X P(\vx) d\mu_X = \int_X P(\vx) d\sigma_X$$
for any polynomial of $P$ of degree $s$; in other words $\mu_X$ and $\sigma_X$
are \emph{$s$-cubature equivalent} \cite{Kuperberg:cubature}.  If we
assume the natural group structure on $T$, then it acts on $X$
in each of the three cases with $Y$ as the set of orbits.  Then
$$\int_X P(\vx) d\sigma_X = \int_X P_T(\vx) d\sigma_X
\quad \int_X P(\vx) d\mu_X = \int_X P_T(\vx) d\mu_X,$$
where $P_T$ is the average of $P$ with respect to the action of $T$.
The polynomial $P_T$ then descends to a polynomial $P_Y$ on $Y$
of degree $t$, and
$$\int_X P_T(\vx) d\sigma_X = P(F_Y)
\qquad \int_X P_T(\vx) d\mu_X = \int_Y P_Y(\vy) d\mu_Y$$
because $\alpha$ preserves measure.

The measure $\sigma_X$ evidently has a twisted product $s$-cubature formula
$F_X = F_T \ltimes F_Y$ given by replacing each fiber by a copy of
$F_T$. (A singular fiber corresponding to a boundary point of $T$ can be
replaced by a projection of $F_T$.) Since $\mu_X$ and $\sigma_X$ are
$s$-cubature equivalent, $F_X$ is a cubature formula for $\mu_X$ as well.
\end{proof}

\begin{remark} The proof of \thm{th:fiber} is analogous
to Sobolev's theorem with the finite group $G$ replaced by the torus
$T$.  Indeed the argument works for any compact group.
\end{remark}

The simplest case of \thm{th:fiber} is the Hopf map $h$. In this case the
theorem says that a $t$-cubature formula $F$ on $\CP^{n-1}$ lifts to a
$(2t+1)$-cubature formula $F'$ on $S^{2n-1}$ with $(2t+2)|F|$ points.  This
relation was also observed by K\"onig \cite{Konig:cubature}.

\begin{corollary} The $n$-sphere $S^n$ has a 7-cubature formula with $O(n^4)$
points for all $n$, more precisely $4n^4(1+o(1))$ points.  The 3-sphere $S^3$
has a $(2s+1)$-cubature formula with 
$$|F| = \begin{cases} (s+1)(s^2+3) & $s$ \mathrm{\ odd} \\
    (s+1)(s^2+s+2) & $s$ \mathrm{\ even} \end{cases}$$
points.
\label{c:fiber} \end{corollary}

\begin{proof} The simplex $\Delta_n$ has a 3-cubature formula with $O(n)$
points \cite{Kuperberg:cubature} constructed using Hadamard designs.  This can
be combined with the 7-design on $T(\SO(2n))$ with $O(n^3)$ points provided by
\thm{th:craig}, for a total of $O(n^4)$ points.  More precisely, the formula on
$\Delta_n$ has points at the corners, each of which lifts to $O(n)$ points; a
point in the center, which lifts to $O(n^3)$ points; and $2n+o(n)$ points on
$\floor{n/2}$-dimensional faces, each of which lift to $2n^3(1+o(1))$ points. 
Only the last family of points is significant and it comprises $4n^4(1+o(1))$
points.

Noskov's formulas from \sec{s:torus} include a $(2s+1)$-design on the square
torus $T(\SO(4))$ with $2(s+1)^2$ points.  When $s$ is odd, this can be
combined with the Gauss-Lobatto $s$-quadrature formula on the interval
$\Delta_2$ with $\frac{s+3}2$ points.  Two of the fibers are circles and
can be replaced by $2(s+1)$ points instead of $2(s+1)^2$ points.  The total is
then $(s+1)(s^2+3)$ points.  When $s$ is even, it can be combined
with the Gauss-Radau $s$-quadrature formula with $\frac{s+2}{2}$ points.
In this case one fiber is a circle.
\end{proof}

\begin{remark} The first part of \cor{c:fiber} actually yields a $7$-design on
$S^{n-1}$ with $O(n^6)$ points whenever there is a Hadamard matrix of order
$n$.  In this case the weights of the $3$-cubature formulas on $\Delta_{n-1}$
are $\frac{2}{n(n+1)(n+2)}$ at the corners, $\frac{n}{2(n+1)(n+2)}$ at the
faces, and $\frac{4n}{(n+1)(n+2)}$ at the center.  Thus the weights are all
commensurable up to a factor of $2n^2$ (note that $n$ is even) and the cubature
formula can be interpreted as a $7$-design with this multiplicity factor. 
Moreover, copies of the lattice formulas on the torus fibers can be shifted to
make the design multiplicity-free.  Better yet, the design need only have
$O(n^4)$ points if, for example $n = 4 \cdot 7^k$.  In this case the prime $p$
used in the proof of \thm{th:craig} can be replaced by the prime power
$7^{k+1}$.  The number of points on each fiber then compensates for all
but a bounded part of the factor of $2n^2$ in the weights.

The previous best construction of $7$-designs on $S^{n-1}$ is due to
Sidelnikov \cite{Sidelnikov:7des} and requires $O(2^{k(k+1)/2})$ points when
$n = 2^k$. \end{remark}

A useful variant of \thm{th:fiber} involves the moment map $\tau_2:\R^{2n} \to
\R^n$ defined by the same formula as $\tau_2:S^{2n-1} \to \Delta_{n-1}$,
namely:
$$\tau_2(x_1,\ldots,x_{2n}) =
    (x_1^2 + x_2^2,x_3^2 + x_4^2,\ldots,x_{2n-1}^2 + x_{2n}^2).$$
This $\tau_2$ takes uniform measure on the ball $B_n$ to uniform
measure on the simplex
$$\Delta'_n = \{\vx | x_k \ge 0, \sum x_k \le 1\}.$$
When $n=2$, Noskov's formulas together with some ad hoc cubature formulas for
the triangle  yield some economical formulas for the 4-ball $B_4$.  For
example, there is a PB 3-cubature formula on the triangle $x,y \ge 0, x+y \le
1$ with points and weights generated from
\begin{align*}
\vp_1 &= (\frac25,\frac25) & w_1 &= \frac{25}{48} \\
\vp_2 &= (\frac{161+17\sqrt{14}}{1344},0) & w_2 &= \frac{16-2\sqrt{14}}{25}
\end{align*}
by switching the coordinates and negating $\sqrt{14}$.  This formula lifts to 1
generic fiber in $B_4$ which can be replaced with 32 points and 4 singular
fibers which are circles and can be replaced with 8 points each. The result is
a PI 7-cubature formula on $B_4$ with 64 points.

Wandzura and Xiao \cite{WX:triangle} found competitive PI $s$-cubature formulas
for $s$ up to 30; \fig{f:xiao} shows one example.  Most of these yield
competitive PI $(2s+1)$-cubature formulas on $B_4$ and $S^5$.  The formulas
could probably be improved further with a search on the triangle that favors
nodes on the edges.

The map $\tau_2:\R^{2n} \to \R^n$ also takes Gaussian measure on $\R^{2n}$ to
exponential measure on $\R_+^n$.  For example, there is a PB exponential
$4$-cubature formula on $\R_+^2$ with points and weights generated from
\begin{align*}
\vp_1 &\approx (1.50766353,1.50766353) & w_1 &\approx 0.354104443 \\
\vp_2 &\approx (6.29508677,1.76717584) & w_2 &\approx 0.00876905581 \\
\vp_3 &\approx (0.285606152,0)         & w_3 &\approx 0.556110610 \\
\vp_4 &\approx (3.27491992,0)          & w_4 &\approx 0.0722468398
\end{align*}
by switching the coordinates.  It lifts to 3 generic fibers with 50 points each
and 4 singular fibers with 10 points each.  The result is a positive Gaussian
7-cubature formula on $\R^4$ with 190 points.

\section{An algebraic lower bound}
\label{s:algebraic}

Let $X$ be the Zariski closure of the support of a measure $\mu$ on $\R^n$ and
let $A$ be the ring of polynomial functions on $X$.  (Recall that the
\emph{Zariski closure} of a set, or closure in the Zariski topology, is the
smallest algebraic variety containing it.)  In other words, $A$ is the quotient
of $\R[\vx]$ by the ideal $I_X$ of polynomials that vanish on $X$.  The ring
$A$ has a degree filtration coming from the degree filtration of $\R[\vx]$. 
Stroud \cite{Stroud:more,Stroud:calc} established an important lower bound on
an arbitrary $2t$-cubature formula $F$ for $\mu$ (not necessarily positive or
interior):

\begin{theorem}[Stroud] If $F$ is a $2t$-cubature formula for $\mu$, then
$$|F| \ge \dim A_{\le t}.$$
\label{th:stroud} \end{theorem}

Mysovskikh \cite{Mysovskikh:exact} observed that this applies to trigonometric
cubature by taking $X = T(\SO(2n))$. (And according to M\"oller
\cite{Moller:lower}, the bound was noted independently in special cases by
other authors, e.g., Radon.)

\begin{proof} Define a bilinear form
$$b:A_{\le t} \times A_{\le t} \to \R$$
by
$$b(P,Q) = \int_X P(\vx)Q(\vx) d\mu.$$
The form $b$ is positive-definite because the integrand of $b(P,P)$ is
non-negative; moreover if the integrand vanishes on $X$, then $P = 0$ as an
element of $A$.  Therefore $b$ is non-degenerate, and its rank is $\dim
A_{\le t}$. On the other hand, the integrand lies in $A_{\le 2t}$, so a
$2t$-cubature formula $F$ leads to the formula
$$b(P,Q) = \sum w_k P(\vp_k)Q(\vp_k).$$
This formula realizes $b$ as a sum of $|F|$ rank 1 forms. Therefore $|F|$
is at least the rank of $b$, as desired.
\end{proof}

An interesting scholium of the proof of \thm{th:stroud} is that if $F$ is a
$2t$-cubature formula, then its points suffice to interpolate polynomials on
$X$ of degree $t$.

It is curiously difficult to improve the Stroud bound for odd-degree cubature.
However, the inference that lower bounds improve mainly in even  degrees is not
consistent with the Hopf fibrations
$$h:S^{2n+1} \to \CP^n \qquad h:T(\SO(2n+2)) \to T(\PSU(n+1)).$$
On the one hand, these maps are quadratic and double the degree of cubature in
passing from the target to the domain; in particular, they do not preserve odd
and even. On the other hand, Sections~\ref{s:projection} and \ref{s:fibration}
together show that cubature in the domain and target are comparably difficult
when $n \gg t$.

The Hopf fibration example suggests a generalization of Stroud's theorem
involving group actions and degree doubling.

\begin{theorem} Let $\mu$ be measure on $\R^n$, and let $X$ be the Zariski
closure of its support.   Let $Y \subset \R^k$ be another affine real algebraic
variety on which a compact group $G$ acts. Let $A$ and $B$ be the rings of
complex-valued polynomials on $X$ and $Y$, and suppose that there is a ring
isomorphism
$$A \stackrel{\cong}{\longto} \Inv_G(B)$$
that doubles the filtration degree of $A$. Let $V$ be a unitary representation
of $G$ and define the filtered vector space
$$M = \Inv_G(B \tensor V).$$
If $F$ is a $t$-cubature formula for $\mu$, then
$$|F| \ge \frac{\dim M_{\le t}}{\dim V}.$$
\label{th:bundle} \end{theorem}

We will always take $Y$ to be the coordinate ring of another algebraic variety
$Y$ which is a principal $G$-bundle over $X$, such that the bundle projection
$\alpha:Y \to X$ is quadratic. The $A$-module $M$ can then be understood as the
space of polynomial sections of a vector bundle $E$ over $X$ with fiber $V$. 
The sections in $M_{\le t}$ then behave like polynomials elements of $A$,
except that their degrees are half-integers.  If $Y = X$ and $G$ is trivial,
then $E$ is the trivial line bundle and \thm{th:bundle} reduces to
\thm{th:stroud}.  The hypotheses of \thm{th:bundle} have been chosen so that
the proof of \thm{th:stroud} generalizes to the case when $E$ is not trivial.

\begin{proof} The vector space $M$ (which is naturally an $A$-module) has an
$A$-valued Hermitian inner product $a$ induced by the Hermitian inner product
on $V$.  More precisely, let $\bV$ be the representation conjugate to $V$ and
let
$$\bM = \Inv_G(\bV \tensor B)$$
be the corresponding conjugate of $M$.  (Note that $A$ and $B$ are both
self-conjugate by hypothesis.)  Let
$$\eps:V \tensor \bV \longto \C$$
be the linearization of the standard Hermitian inner product on $V$, and let
$$m:B \tensor B \longto B$$
be the linearization of multiplication on $B$. Let $a'$ be the composition
$$B \tensor V \tensor \bV \tensor B \stackrel{I \tensor \eps \tensor I}
    {\longto} B \tensor B \stackrel{m}{\longto} B.$$
We can restrict the domain to 
$$M \tensor \bM = \Inv_G(B \tensor V) \tensor \Inv_G(\bV \tensor B).$$
Since the restricted domain is $G$-invariant, we can then restrict the target
to $A$.  Let $a$ be this restriction of $a'$.  Although given as a linear map
on $M \tensor \bM$, it can be reinterpreted as a Hermitian inner product on
$M$.  In more geometric terms, if $M$ comes from a bundle $E$ over $X$ with
fiber $V$, then $a(f,g)$ is the pointwise inner product of two sections $f$ and
$g$ of $E$.

Note that $a$ is positive-definite in the sense that
$$a(f,f)(\vx) \ge 0$$
for all $\vx \in X$, and if $a(f,f) = 0$, then $f = 0 \in M$.
The rest of the proof follows that of \thm{th:stroud}:
Define a complex-valued Hermitian inner product $b$ on $M_{\le t}$ by
$$b(f,g) = \int_X a(f(\vx),g(\vx)) d\mu.$$
Then $b$ is also positive-definite, because $a$ is positive-definite and $\mu$
is Zariski-dense in $X$.  Thus, $b$ has rank $\dim M_{\le t}$. A cubature
formula $F$ realizes $b$ as a sum of $|F|$ terms of rank at most $\dim V$.
\end{proof}

We state three special cases of \thm{th:bundle} as corollaries:

\begin{corollary} If $|F|$ is a $(2t+1)$-cubature
formula on $\CP^n$, then
$$|F| \le \binom{n+t}{n}\binom{n+t+1}{n}.$$
\end{corollary}
\begin{proof} Let $Y$ be $S^{2n+1}$ and let $G$ to $S^1 \subset \C$ acting by
complex rotation on $\C^{n+1} \supset S^{2n+1}$. The bundle projection
$\alpha:Y \to X$ is the Hopf map $h$.  Let $V = L_1$ be the tautological
representation of $S^1$, so that $E$ is the tautological line bundle on
$\CP^n$.  The space $M_{\le 2t+1}$ is explicitly realized as the space of
homogeneous polynomials in $\vz$ and $\bar{\vz}$ of bidegree $(t+1,t)$.
The result follows by noting that
$$\dim M_{\le 2t+1} = \binom{n+t}{n}\binom{n+t+1}{n}$$
and that $\dim L_1 = 1$.
\end{proof}

\begin{corollary} If $|F|$ is a $t$-cubature formula on $T(\PSU(n+1))$, then
$$|F| \ge \frac{n^t(1-o(1))}{\ceil{t/2}!\floor{t/2}!}.$$
\label{c:molpsu}
\end{corollary}

\begin{proof} Let $Y$ be the torus $T(\SO(2n+2))$ and let $G = S^1$ again act by
complex rotation in $\C^{n+1}$. In this case $M_{\le 2t+1}$ is spanned by the
space of monomials in $\vz$ and $\bar{\vz}$ of bidegree $(s+1,s)$ with $s \le
t$ and with the relation
$$z_k\bz_k = 1$$
for all $k$. Its dimension is the number of points in the Minkowski difference
$$\Delta_n^{(t+1)} - \Delta_n^{(t)},$$
where $\Delta_n^{(t)}$ is the discrete simplex defined in the proof of
\thm{th:craig}.  This is very similar to \thm{th:stroud} for $2t$-cubature,
because
$$\dim A_{\le t} = |\Delta_n^{(t)} - \Delta_n^{(t)}|.$$
There is no concise formula for either number, but there is a concise estimate
for fixed $t$ in the limit $n \to \infty$.  If $E$ is either the trivial bundle
when $t$ is even or the bundle $L_1$ (restricted from $\CP^n$) when $t$ is odd,
then
$$\dim M_{\le t} \approx \binom{n+1}{\ceil{t/2},\floor{t/2},n+1-t} \approx
\frac{n^t}{\ceil{t/2}!\floor{t/2}!}$$
as $t \to \infty$, as desired.
\end{proof}

\begin{remark} When $F$ is a lattice formula, \cor{c:molpsu} is equivalent
to Minkowski's classic upper bound on the density $\hF$ as a lattice packing
of the discrete simplex $\Delta^{(t)}_n$.  This and the fact that the Hopf
fibration is quadratic led the author to \thm{th:bundle}.
\end{remark}

\begin{corollary} If $\mu$ is a Zariski-dense measure on $S^n$
and $F$ is a $(2t+1)$-cubature formula for $\mu$, then
$$|F| \ge 2\binom{n-1+t}t.$$
\label{c:molsphere} \end{corollary}

\begin{proof} The idea is to let $Y = \Spin(n+1)$ and $G = \Spin(n)$, where
$\Spin(n)$ is the connected Lie group that double covers $\SO(n)$.  Then
$$X = \Spin(n+1)/\Spin(n) = \SO(n+1)/\SO(n) = S^n.$$
Our choice for the representation $V$ is the spinor representation of
$\Spin(n)$ when $n$ is odd and a semispinor representation when $n$ is even. 
The rest of the argument is a review of standard but non-trivial representation
theory \cite{Varadarajan:gtm,Koike:spinor}:

\begin{description}
\item[1.] The irreducible (unitary) representations of the spin groups
$\Spin(2n)$ and $\Spin(2n+1)$ are both indexed by vectors
$$\vlam = a_1\vlam_1 + a_2\vlam_2 + \ldots + a_n\vlam_n$$
with each $a_k$ a non-negative integer.  The representation with highest weight
$\vlam$ can be called $V_G(\vlam)$ with $G = \Spin(2n)$ or $G =
\Spin(2n+1)$.

\item[2.] The representation $V_G(t\vlam_1)$
is realized as the space $A_t$ of homogeneous polynomials of degree
$t$ on $S^{2n-1}$ or $S^{2n}$.  The representation $V_{\Spin(2n+1)}(\vlam_n)$
is the spinor representation of dimension $2^n$.  The representations
$V_{\Spin(2n)}(\vlam_{n-1})$ and $V_{\Spin(2n)}(\vlam_n)$
are the semispinor representations of dimension $2^{n-1}$.

\item[3.] For any choice of fiber $V = V_{\Spin(n)}(\vlam)$, the $A$-module $M$
is the induced representation
$$M = M(n,\vlam) = \Ind_{\Spin(n)}^{\Spin(n+1)} V_{\Spin(n)}(\vlam).$$

\item[4.] The characters of $V(\vlam)$ for all $\vlam$ generate the natural
real algebraic structure on $\Spin(2n)$ or $\Spin(2n+1)$ (indeed on any compact
Lie group).  One degree filtration is defined by letting the characters of every
$V(\vlam_k)$ have degree 2, except for the spinor or semispinor
representations, which have degree 1.

\item[5.] The structure of $M(n,\vlam)$ can be computed by branching formulas
for the restriction of an irrep of $\Spin(n+1)$ to $\Spin(n)$, together with
induction-restriction duality. These restriction formulas were computed by
Koike \cite[Thms.~11.2\&11.3]{Koike:spinor}.  When $\vlam = \vlam_n$, Koike's
formulas together with the degree filtrations yield
$$M(2n+1;\vlam_n)_{\le 2t+1} \cong \bigoplus_{s \le t}
    V(2n+1;s\vlam_1 + \vlam_n)$$
and
$$M(2n;\vlam_n)_{\le 2t+1} \cong \bigoplus_{s \le t}
    V(2n;s\vlam_1 + \vlam_{n-1}) \oplus V(2n;s\vlam_1 + \vlam_n).$$

\item[6.] Finally by the Weyl dimension formula,
$$\dim V(2n+1;s\vlam_1 + \vlam_n) = 2^n\binom{2n-1+s}{s}$$
and
$$\dim V(2n;s\vlam_1 + \vlam_n) = \dim V(2n;s\vlam_1 + \vlam_n)
    = 2^{n-1}\binom{2n-2+s}{s}.$$
\end{description}

Combining the dimension formulas yields
$$|F| \ge \frac{\dim M(n;\vlam)_{\le 2t+1}}{\dim V}
    = \sum_{s \le t} 2\binom{n-2+s}{s} = 2\binom{n-1+t}{t},$$
as desired.
\end{proof}

\begin{remark} \cor{c:molsphere} matches the M\"oller bound \cite{Moller:lower} for cubature
on $S^n$, but it is more general because the measure $\mu$ need not be
centrally symmetric.
\end{remark}

\begin{fullfigure*}{f:xiao}
    {A 175-point 30-cubature formula found by Wandzura and Xiao \cite{WX:triangle}}
\includegraphics{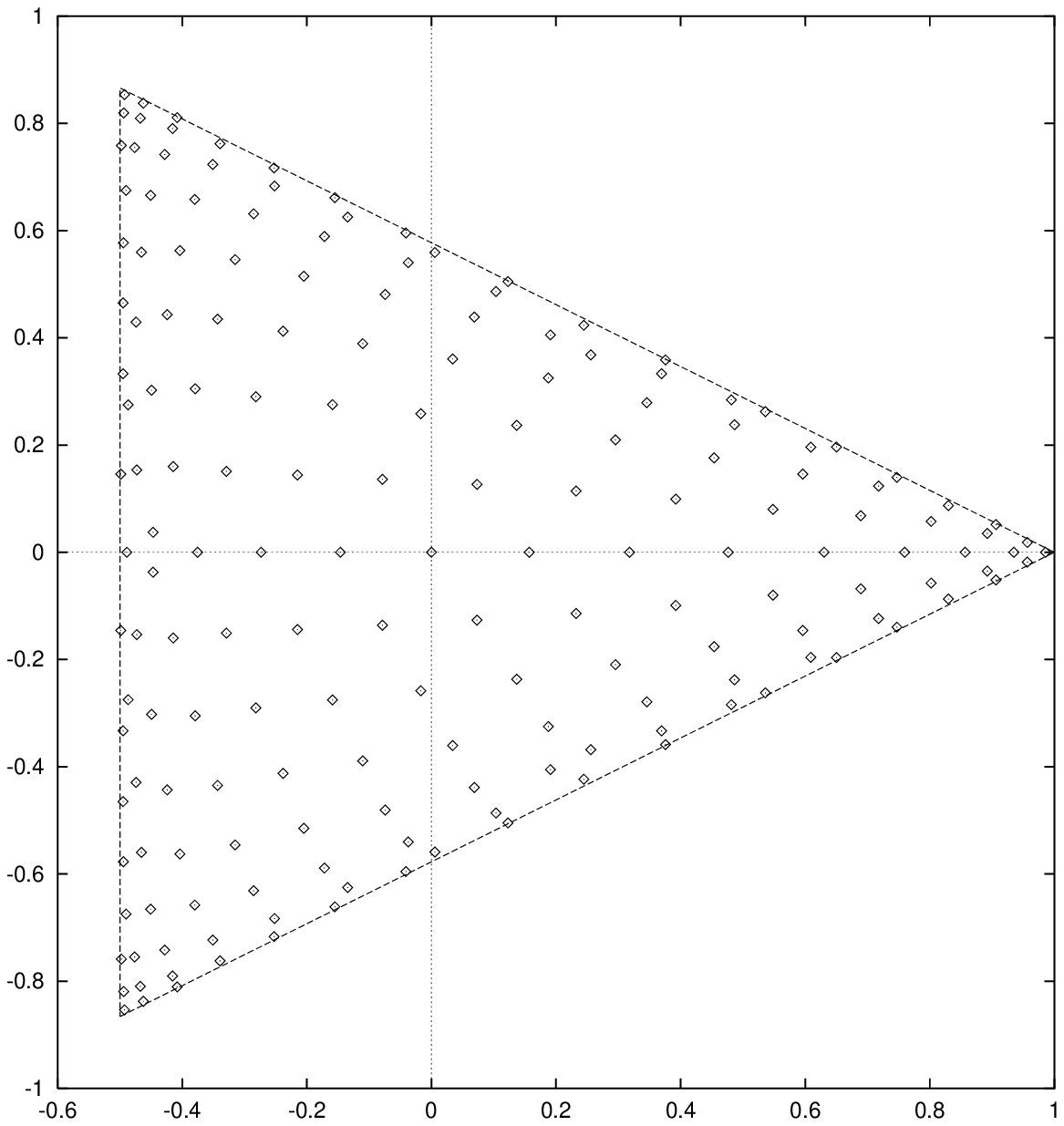}
\end{fullfigure*}

\section{A local lower bound}
\label{s:local}

Our final application of moment maps is to help establish a local lower bound
on the density of points of a PI or PB cubature formula on the simplex
$\Delta_n$.  The bound was originally inspired by PI cubature formulas due to
Wandzura and Xiao \cite{WX:triangle} which were found by simulated annealing. 
As in the example shown in \fig{f:xiao}, the points in these formulas
accumulate transversely at the edges of the triangle. Another related result is
that the limiting density of the points of Gauss-Legendre quadrature (\ie, the
zeros of Legendre polynomials) is $\frac{1}{\pi\sqrt{1-x^2}}$
\cite{Dehesa:transport}.  This density can be interpreted as the linear
projection of uniform measure on a circle, which is related to Archimedes'
moment map (\fig{f:xu}).

\thm{th:local} and \cor{c:limit} establish a lower bound on the limiting
density of any sequence of PI and PB formulas on $\Delta_n$ that generalizes
the limiting density of Legendre zeros. Moreover, if the local density is high
in certain regions, in particular near the vertices of $\Delta_n$, then the
weights there must be low.  By this reasoning, \thm{th:sharp} and
\cor{c:design} establish that a $t$-design on $\Delta_n$ requires many more
points than an efficient $t$-cubature formula does) as $t \to \infty$ (namely
$O(t^{2n})$ points versus $O(t^n)$ points).  Along the way, \thm{th:sharp}
establishes that Gaussian quadrature for an arbitrary weight function is very
sharply locally optimal among all positive quadrature formulas. Finally
\sch{s:simple} generalizes the results for uniform measure on $\Delta_n$ to
uniform measure on an arbitrary simple convex polytope.

\begin{theorem} A PI or PB $2t-1$-cubature formula $F$ on the simplex
$\Delta_n$ is an $\eps$-net with respect to the metric
$$ds^2 = \frac{dx_0^2}{2x_0} + \frac{dx_1^2}{2x_1} + \ldots
    + \frac{dx_n^2}{2x_n}$$
in barycentric coordinates, where $\cos 2\eps$ is the highest zero of the
Jacobi polynomial $P^{(n-1,0)}_t(x)$.
\label{th:local} \end{theorem}

In the proof and later, we will abbreviate $(n-1,0)$ as ``$\#$''in
superscripts.

\begin{fullfigure}{f:island}
    {A polynomial on $\Delta_n$ with a small positive island}
\pspicture(-2,-1.25)(2,2.25)
\pspolygon(0,2)(-1.732,-1)(1.732,-1)
\rput{-50}(.4,.7){\psellipse(0,0)(.4,.175)}\rput(.4,.7){$+$}
\rput(1.04,-.6){$-$}\rput(-.3,.6){$-$}\rput(-.7,-.5){$-$}\rput(0,0){$-$}
\endpspicture
\end{fullfigure}

\begin{proof} The idea of the proof is to find, for each $\vp \in \Delta_n$ and
each $\eps' > \eps$, a $P$ of degree $2t-1$ on $\Delta_n$ such that 
$$\int_{\Delta_n} P(\vx) d\vx > 0,$$
but $P(\vx) > 0$ only when $\vx \in \Delta_n$ is in the $\eps'$-ball
$B_{\eps'}(\vp)$ around $\vp$.  We can call this ball the \emph{positive
island} of $P(\vx)$; see \fig{f:island}.  The existence of such a polynomial
$P$ forces $F$ to have an evaluation point in $B_\eps(\vp)$, for otherwise
$P(F) \le 0$.

We first claim that the stated metric is the distance between fibers of the
moment map $\pi$ with respect to the Fubini-Study metric on $\CP^n$.  To see
this it suffices to check the following: The real locus $\RP^n \subset \CP^n$
is perpendicular to the fibers of $\pi$ and meets each generic fiber $2^n$
times.  Indeed, $\pi$ is a bijection on the orthant $\RP_{\ge 0}^n$ with
non-negative projective coordinates.  Moreover, $\RP_{\ge 0}^n$ is isometric to
the orthant $S^n_{\ge 0}$ of a unit $n$-sphere, and the restriction of $\pi$
agrees with the restriction of the Xu map $\tau_1$.  The metric $ds^2$ on
$\Delta_n$ is exactly the push-forward of the standard metric on $S^n_{\ge 0}$
under $\tau_1$.  See \fig{f:xu} for an example.

\begin{fullfigure}{f:xu}{The moment map $\pi$ restricts to the Xu map $\tau_1$}
\pspicture(-2.3,-2.1)(4.85,2.1)
\pscircle(0,0){2}
\psbezier[linestyle=dashed](-0.518,-1.783)(-0.518,-1.893)(-0.286,-1.982)(0.000,-1.982)
\psbezier[linestyle=dashed](0.518,-1.783)(0.518,-1.893)(0.286,-1.982)(0.000,-1.982)
\psbezier[linestyle=dashed](-0.518,-1.783)(-0.518,-1.673)(-0.286,-1.584)(0.000,-1.584)
\psbezier[linestyle=dashed](0.518,-1.783)(0.518,-1.673)(0.286,-1.584)(0.000,-1.584)
\psbezier(-1.286,-1.414)(-1.286,-1.687)(-0.710,-1.909)(0.000,-1.909)
\psbezier(1.286,-1.414)(1.286,-1.687)(0.710,-1.909)(0.000,-1.909)
\psbezier[linestyle=dashed](-1.286,-1.414)(-1.286,-1.141)(-0.710,-0.920)(0.000,-0.920)
\psbezier[linestyle=dashed](1.286,-1.414)(1.286,-1.141)(0.710,-0.920)(0.000,-0.920)
\psbezier(-1.813,-0.780)(-1.813,-1.165)(-1.001,-1.477)(0.000,-1.477)
\psbezier(1.813,-0.780)(1.813,-1.165)(1.001,-1.477)(0.000,-1.477)
\psbezier[linestyle=dashed](-1.813,-0.780)(-1.813,-0.395)(-1.001,-0.083)(0.000,-0.083)
\psbezier[linestyle=dashed](1.813,-0.780)(1.813,-0.395)(1.001,-0.083)(0.000,-0.083)
\psbezier(-2.000,0.000)(-2.000,-0.425)(-1.105,-0.769)(0.000,-0.769)
\psbezier(2.000,0.000)(2.000,-0.425)(1.105,-0.769)(0.000,-0.769)
\psbezier[linestyle=dashed](-2.000,0.000)(-2.000,0.425)(-1.105,0.769)(0.000,0.769)
\psbezier[linestyle=dashed](2.000,0.000)(2.000,0.425)(1.105,0.769)(0.000,0.769)
\psbezier(-1.813,0.780)(-1.813,0.395)(-1.001,0.083)(0.000,0.083)
\psbezier(1.813,0.780)(1.813,0.395)(1.001,0.083)(0.000,0.083)
\psbezier[linestyle=dashed](-1.813,0.780)(-1.813,1.165)(-1.001,1.477)(0.000,1.477)
\psbezier[linestyle=dashed](1.813,0.780)(1.813,1.165)(1.001,1.477)(0.000,1.477)
\psbezier(-1.286,1.414)(-1.286,1.141)(-0.710,0.920)(0.000,0.920)
\psbezier(1.286,1.414)(1.286,1.141)(0.710,0.920)(0.000,0.920)
\psbezier[linestyle=dashed](-1.286,1.414)(-1.286,1.687)(-0.710,1.909)(0.000,1.909)
\psbezier[linestyle=dashed](1.286,1.414)(1.286,1.687)(0.710,1.909)(0.000,1.909)
\psbezier(-0.518,1.783)(-0.518,1.673)(-0.286,1.584)(0.000,1.584)
\psbezier(0.518,1.783)(0.518,1.673)(0.286,1.584)(0.000,1.584)
\psbezier(-0.518,1.783)(-0.518,1.893)(-0.286,1.982)(0.000,1.982)
\psbezier(0.518,1.783)(0.518,1.893)(0.286,1.982)(0.000,1.982)
\psbezier[linewidth=1pt](1.600,-0.462)(1.600,0.558)(0.884,1.591)(0.000,1.846)
\psbezier[linewidth=1pt](1.600,-0.462)(1.600,-1.481)(0.884,-2.101)(0.000,-1.846)
\qdisk(0.000,1.846){.1}\qdisk(0.000,-1.846){.1}
\psline{<-}(1.369,0.858)(1.969,1.458)\rput[bl](2.019,1.458){$\RP_{\ge 0}^1$}
\rput[br](-1.7,1.2){$\CP^1$}
\psline{->}(2.5,0)(4,0)\rput[b](3.25,.13){$\pi$}
\psline(4.5,-2)(4.5,2) \psline(4.4,-2)(4.6,-2) \psline(4.4,2)(4.6,2)
\endpspicture\eatline
\end{fullfigure}

Consider the linear projection $\alpha:\Delta_n \to [-1,1]$ given by
\eq{e:alpha}{\alpha(\vx) = 2x_0-1.}
The map $\alpha$ sends uniform measure on $\CP^n$ to the measure
$$\mu(x) = n2^{1-n}(1-x)^{n-1}.$$
The $t$th orthogonal polynomial with respect to this measure $\mu$ is the
Jacobi polynomial $P^\#_t = P^{(n-1,0)}_t$. Let $p^\#_t$ be its highest zero.

Define a polynomial $Q_\delta:\CP^n \to \R$ by
\eq{e:qdef}{Q_\delta(\vz) = Q_\delta(x) =
    \frac{P^\#_t(x)^2(x-p^\#_t+\delta)}{(x-p^\#_t)^2},}
where $x = \alpha(\pi(\vz))$ and $\delta > 0$.  It has degree $2t-1$ as a
polynomial in $x$, as well as a polynomial on $\CP^n$. Moreover, $Q_\delta$
vanishes at the zeros of $P^\#_t$, except at the highest zero, at which its
value is positive.  Therefore by Gaussian quadrature (!) with respect to the
measure $\mu$,
$$\int_{\CP^n} Q_\delta(\vz) d\vz = \int_{-1}^1 Q_\delta(x) d\mu > 0.$$
At the same time, $Q$ is non-positive outside of the region
$$x > p^\#_t - \delta.$$
This region corresponds to the ball of radius $\eps'$ around
$(1:0:0:\cdots:0)$, with
$$2(\cos \eps')^2 - 1 = \cos 2\eps' = p^\#_t - \delta.$$
This can be confirmed by comparing with the orthant $\RP_{\ge 0}^n$ from
mentioned previously.  Note that $\eps' \to \eps$ as $\delta \to 0$.

Given $\vq \in \CP^n$, define $Q_{\delta,\vq}$ by rotating $Q_{\delta}$ by some
isometry of $\CP^n$ that takes $(1:0:0:\cdots:0)$ to $\vq$. Define
$Q^T_{\delta,\vp}:\Delta_n \to \R$, where $\vp = \pi(\vq)$, by averaging $Q$ over torus fibers:
$$Q^T_{\delta,\vp}(\vx) = \frac1{|\pi^{-1}(\vx)|}
    \int_{\pi^{-1}(\vx)} Q_{\delta,\vp}(\vz) d\vz.$$
Then
$$\int_{\Delta_n} Q^T_{\delta,\vp}(\vx) d\vx =
    \int_{\CP^n} Q_{\delta,\vp}(\vz) d\vz > 0,$$
and $Q^T_{\delta,\vp}$ is non-positive outside of the ball of radius $\eps'$
around $\vp = \pi(\vq)$ in the induced metric on $\Delta_n$. Thus,
$Q^T_{\delta,\vp}$ has the desired properties.
\end{proof}

\begin{remark} A somewhat weaker version of \thm{th:local} holds when $F$ is
positive and exterior, but with real evaluation points.  Polynomials similar to
$Q^T_\delta$ can be constructed directly as products of factors that vanish on
quadratic surfaces in $\R^n \supset \Delta_n$, with only one unsquared factor
that vanishes on the boundary of $B_{\eps'}(\vp)$.  As it happens, the boundary
of $B_{\eps'}(\vp)$ is a quadratic surface.  We did not refine this
sketched argument into a proof with explicit estimates.
\end{remark}

\begin{corollary} Any sequence of PI or PB $t$-cubature formulas on $\Delta_n$
has limiting point density $\Omega(t^n\prod_k x_k^{-1/2})$, where $\vx \in
\Delta_n$ is fixed and given in barycentric coordinates, and $t \to \infty$. 
\label{c:limit}\end{corollary}

\begin{proof} The corollary follows from computing the volume form
corresponding to the metric $ds^2$ in the statement of \thm{th:local} and
estimating the covering radius $\eps$.   The asymptotic behavior of zeros of
Jacobi polynomials is given in Abramowitz and Stegun \cite[p.
787]{AS:handbook}.  The key step in the estimate is the limit
\eq{e:bessel}{\lim_{t \to \infty} \frac{P^{(a,b)}_t(\cos \frac{\theta}{t})}
    {P^{(a,b)}_t(1)} = 2^a\theta^{-a}a!J_a(\theta),}
where $J_a(z)$ is the ordinary Bessel function of the first kind.  Convergence
to the limit is analytic in $\theta$.  Thus 
$$\lim_{t \to \infty} t\theta^{(a,b)}_{t,t+1-k} = j_{a,k}$$
for every fixed $k$, where $\cos \theta^{(a,b)}_{t,k}/t$ is the $k$th zero of
$P^{(a,b)}_t(x)$ and $j_{a,k}$ is the $k$th zero of $J_a(x)$.  The estimate can
be established directly in our geometry by noting that
$P^\#_t(2|z_0|^2-1)$
is a harmonic function on $\CP^n$.  The harmonic equation on $\CP^n$ is then
locally approximated by the harmonic (or Helmholtz) equation on $\R^{2n}$,
whose radial solutions are derived from Bessel functions.

In our case,
$$\cos 2\eps = \cos \frac{\theta^\#_{t,t}}{t},$$
for $(2t-1)$ cubature.  So
$$\eps = \frac{j_{n-1,1}}{2t}(1+o(1)) = \Theta(t^{-1}),$$
which is also $\Theta(t^{-1})$ for $t$-cubature.
\end{proof}

\begin{theorem} Let $\mu$ be an arbitrary normalized measure on $\R$ whose
support has at least $2t$ points.  Let $p_1,\ldots,p_t$ and $w_1,\ldots,w_t$ be
the points and weights of Gaussian $t$-quadrature for the measure $\mu$.  Let
$F$ be a positive $t$-quadrature formula for $\mu$. Then for each $1 \le k
\le t$, $F$ has at least one point in the half-open interval $(p_{k-1},p_k]$,
where $p_0 = -\infty$.  Moreover, the total weight of all points in
$(-\infty,p_1]$ is at most $w_1$, with equality if and only if $F$ is the
Gaussian quadrature formula.
\label{th:sharp} \end{theorem}

Note that \thm{th:sharp} is further sharpened by symmetry: $F$ must also have
at least one point in each half-open interval $[p_k,p_{k+1})$, with $p_{t+1} =
\infty$, and its total weight in $[p_t,\infty)$ is at most $w_t$.

\begin{proof}
Let $\phi_t(x)$ be the $t$th orthonormal polynomial with respect to $\mu$ (with
either sign), and let $A_t$ be the leading coefficient of $\phi_t(x)$. If
$k=1$, let
$$P(x) = \frac{\phi_t(x)^2(p_1+\delta_1-x)}{(x-p_1)^2}$$
with $\delta_1 > 0$. If $k > 1$, let
$$P(x) = \frac{\phi_t(x)^2(x-p_{k-1}-\delta_0)
    (p_k+\delta_1-x)}{(x-p_{k-1})^2(x-p_k)^2}$$
with $\delta_1 \gg \delta_0 > 0$.  In both cases,
$$\int_\R P(x) d\mu > 0$$
by Gaussian quadrature, while $P$ is only positive on the interval
$(p_{k-1}+\delta_0,p_k+\delta_1)$.  Therefore $F$ has at least one point in
this interval.  Since $F$ only has finitely many points, the limit $\delta_1
\to 0$ establishes that $F$ has a point in $(p_{k-1},p_k]$.

For the second claim, let 
$$P(x) = \frac{\phi_t(x)^2}{(x-p_1)^2}.$$
Then by Gaussian quadrature,
$$\int_\R P(x) d\mu = w_1P(p_1).$$
Let $q_1,\ldots,q_k$ be the points of $F$ which are at most $p_1$, and let
$v_1,\ldots,v_k$ be their weights.  Then
$$\int_\R P(x) d\mu = P(F)
    \ge \sum_{j=1}^k w_j P(q_j) \ge P(p_1) \sum_{j=1}^k w_j.$$
The first inequality holds because $P$ is non-negative; the second because $P$
decreases on $(-\infty,p_1]$.
\end{proof}

\begin{corollary} The least weight of any positive $t$-cubature formula on
$\Delta_n$ (with uniform measure) is $O(t^{-2n})$, uniformly in $t$. Any
$t$-design on $\Delta_n$ has $\Omega(t^{2n})$ points.
\label{c:design} \end{corollary}
\begin{proof}
If $F$ is a $t$-cubature formula on $\Delta_n$, the map $\alpha$ (see
equation~\eqref{e:alpha}) sends it to a $t$-quadrature formula $\alpha(F)$
on $[-1,1]$ with Jacobi-polynomial measure. If $F$ is positive, then the
least weight of $\alpha(F)$ is at least that of $F$.  On the other hand,
\thm{th:sharp} establishes that the least weight $\alpha(F)$ is at least the
last Christoffel weight $w_t$.  

The first claim follows by estimating this weight. One of the standard formulas
for the general Christoffel weight $w_k$ is
$$w_k = -\frac{A_{t+1}||\phi_t(x)||_\mu^2}{A_t
    \phi_t'(p_k)\phi_{t+1}(p_k)},$$
where $\phi_t(x)$ is the $t$th orthogonal polynomial, $A_t$ is its leading
coefficient, and $p_k$ is its $k$th root. In our case, $\phi_t = P^\#_t$, $p_k
= p^\#_{t,k}$, and $k=t$.  We compute:
\begin{align*}
||P^\#_t||_\mu^2 &= \frac{n}{2t+n} = \Theta(t^{-1}) \\
    A_t &= 2^{-t}\binom{n-1+2t}{t} = \Theta(2^t).
\end{align*}
To estimate $(P^\#_t)'(p^\#_{t,t})$ and $P^\#_{t+1}(p^\#_{t,t})$, we again
appeal to the limit in equation \eqref{e:bessel}. Differentiating both sides by
$\theta$, we obtain
$$\lim_{t \to \infty} - \frac{(P^\#_t)'(\cos \theta/t)
    (\sin \theta/t)}{tP^\#_t(1)} = -a2^a\theta^{-a}a!J'_a(\theta).$$
Note that
$P^\#_t(1) = \binom{t+n-1}{t} = \Theta(t^{n-1})$.
For a fixed value of $\theta$, the various parts of the limit yield
$$(P^\#_t)'(\cos \frac{\theta}{t}) = \Theta(t^{n+1}).$$
By the same token
$$(P^\#_t)'(\cos \frac{\theta_t}{t}) = \Theta(t^{n+1})$$
when $\theta_t$ approaches a fixed value of $\theta$, as is the case when
$\theta_t = \theta^\#_{t,t}$ is given by
$$p^\#_{t,t} = \cos \frac{\theta^\#_{t,t}}{t}.$$
By a similar calculation,
$$(P^\#_{t+1})(p^\#_{t,t}) = \Theta(t^{n-2}).$$
The conclusion is that
$$w_t = \Theta(t^{2n}),$$
as desired.
\end{proof}

\begin{scholium} Let $K \subset \R^n$ be a convex $n$-polytope with $N$
facets.  Let $F$ be a $t$-cubature formula on $K$ with uniform measure.  If
$F$ is PI or PB and if $K$ is simple, then $F$ is an $\eps$-net with
respect to the metric
$$ds^2 = \frac{dx_0^2}{x_1} + \frac{dx_2^2}{x_2} + \ldots
    + \frac{dx_N^2}{x_N},$$
where $\eps = O(1/t)$.  If $F$ is positive, then its least weight is
$O(t^{-2n})$. If it is a $t$-design, then it has at least $O(t^{2n})$ points.
\label{s:simple} \end{scholium}

\begin{fullfigure}{f:facets}{A simplex $L^{-1}(\Delta_n)$ whose facets
    contain facets of $K$ that meet at $\vx$}
\pspicture(-3.5,-1.6)(3.5,1.9)
\pspolygon(3.236;198)(3.236;342)(1.236;90)
\psline(1.236;18)(1.236;306)
\psline(1.236;162)(1.236;234)
\rput(0,0){$K$}
\rput[t](-1.3,-1.2){$L^{-1}(\Delta_n)$}
\rput[b](0,1.45){$\vx$}
\qdisk(1.236;90){.07}
\endpspicture\eatline
\end{fullfigure}

\begin{proof}(Sketch)
The proof of \thm{th:local} retains its strength if $K \subseteq \Delta_n$ and
we pass from $\Delta_n$ to $K$, provided that the positive island of the
polynomial $Q^T_{\delta,\vp}$ lies within $K$.  In this case
$$\int_K Q^T_{\delta,\vp}(\vx) d\vx
    \ge \int_{\Delta_n} Q^T_{\delta,\vp}(\vx)d\vx > 0.$$
In order to properly position $Q^T_{\delta,\vp}$ for all $\vp \in K$, we need
several embeddings of $K$ into $\Delta_n$.  For each vertex $\vx \in K$, choose
a linear embedding $L$ that sends $\vx$ to some vertex of $\Delta_n$, and that
sends the facets incident to $\vx$ to facets of $\Delta_n$.  (Equivalently, for
each vertex $\vx \in K$, choose a simplex $L^{-1}(\Delta_n) \supseteq K$ whose
facets includes all facets of $K$ that meet at $\vx$.  See \fig{f:facets}.)
Then there exists a finite set of $L$ such that the positive islands of
polynomials of the form $Q^T_{\delta,vq} \circ L$ together cover $K$.  The
formula $F$ must have a point in each island, which establishes that $F$ is an
$\eps$-net.

Similarly, the proofs of \thm{th:sharp} and \cor{c:design} retain their
strength if uniform measure on $K$ projects by a map $\alpha$ to a measure
$\nu$ on $[-1,1]$ which is dominated by
$$\mu(x) = 2^{-n}n(1-x)^{n-1}$$
and which agrees with $\mu$ in a neighborhood of $1$. (Of course $\mu$
cannot dominate $\nu$ if $\nu$ is normalized, so this condition on $\nu$ must
be dropped.)  In this case
$$\int_\R P(x) d\nu \ge \int_\R P(x) d\mu > 0$$
for the first half of \thm{th:sharp} for the interval $(-\infty,p_1)$, while
$$\int_\R P(x) d\nu \le \int_\R P(x) d\mu$$
for the second half of \thm{th:sharp}.  A suitable projection $\alpha$ can be
realized by positioning $K$ in $\Delta_n$ so that it touches the vertex $x_0 =
1$, and then restricting the usual map $\alpha$ to $K$.
\end{proof}


\section{Other comments}
\label{s:other}

In this article we have studied the toric moment map on $\CP^n$, and on $\C^n$
restricted to $S^{2n-1}$ (which can be interpreted as the level surface of an
invariant Hamiltonian on $\C^n$) as it applies to the cubature problem. Many of
the constructions apply equally well to arbitrary toric varieties. To begin
with, every complex projective variety $X$ inherits both a metric and an affine
real structure from $\CP^n$.  If $X$ is toric, it also has a volume-preserving
moment map whose image is a centrally symmetric polytope. However, the variety
$X$ rarely has much more symmetry than its moment map image.

The duality between toric cubature (in particular trigonometric cubature) and
lattice packings explored in \sec{s:torus} suggests a different limit of the
cubature problem.  Let $K \subset \R^n$ be a centrally symmetric convex body. 
For simplicity let $F = \{\vp_a\}$ be a periodic discrete subset of $\R^n$ with
a periodic weight function $\vp_a \mapsto w_a$.  Since it is periodic, it has a
well-defined Fourier transform $\hF$. In this setting, $F$ is a Fourier
$K$-cubature formula if and only if $\hF$ is disjoint from the interior of
$K$. The (continuous) Fourier cubature problem is to minimize the density of
$F$ among all $K$-cubature formulas or all positive $K$-cubature formulas.  If
$F$ is a lattice with equal weights, then $\hF$ has the same property and
Fourier $K$-cubature problem reduces to finding the best lattice packing of
$K$.  It would be interesting to find examples of non-lattice formulas that are
better than the best lattice formula.

We conjecture that a version of \cor{c:molsphere} holds, using \thm{th:bundle}
and the same spinor bundles, for any centrally symmetric subvariety $X \subset
S^n$.  That is, we conjecture M\"oller's bound for these varieties, even when
the measure $\mu$ is not centrally symmetric.

\thm{th:local} shows why some tempting approaches to construct efficient PI or
PB formulas on the simplex $\Delta_n$, even the triangle $\Delta_2$, are bound
to fail.  For example, if the points of a putative cubature formula $F$ are
fixed in advance, the question of whether it admits non-negative weights for
$t$-cubature reduces to linear programming. But if the points are arranged in
some lattice with spacing $1/k$, \thm{th:local} shows that the weights can only
be non-negative if $k = \Omega(t^2)$, so that $|F| = \Omega(t^{2n})$.

We believe that the requirement that $K$ be simple in \sch{s:simple} is not
essential. More generally we conjecture that similar results hold if $K$ is not
convex.  We also conjecture that the bounds in \thm{th:local} and
\sch{s:simple} are sharp to within a constant factor. The cubature formulas
found by Wandzura and Xiao support this conjecture, at least when $K =
\Delta_n$.

The proof of \thm{th:local} was partly inspired by the linear programming
method to bound kissing numbers, $t$-designs, and sphere packings
\cite{Delsarte:linear,DGS:spherical,KL:sphere,OS:unit,CE:new}. Xu observed that
the method for $t$-designs also yields bounds on PI $t$-cubature
\cite{Xu:lower}. In fact it yields an upper bound on the $\ell_2$ norm of the
weights of a PI $t$-cubature formula, which implies a lower bound on the number
of points.  We conjecture that linear programming methods could be used to
improve the constants in \thm{th:local}.

Krylov \cite{Krylov:calc} established that if $\{F_t\}$ is a sequence of
interior $t$-cubature formulas for a measure $\mu$, then $\{f(F_t)\}$ converges
to $\int_X f(\vx) d\mu$ for every continuous $f$ if and only if  the $\ell_1$
norm of the weights of $F_t$ is bounded as $t \to \infty$. We conjecture then
that \thm{th:local} still holds assuming a bound on the $\ell_1$ norm of the
coefficients of $F$ instead of assuming that $F$ is positive.

\acknowledgments

The author would like to thank Yael Karshon, W{\l}odzimierz Kuperberg, Thomas
Strohmer, and Hong Xiao for useful discussions. The author would especially
like to thank Noam Elkies and Eric Rains, who obtained some of the examples and
provided many other useful comments.

\bibliography{mg,me,shared}

\providecommand{\bysame}{\leavevmode\hbox to3em{\hrulefill}\thinspace}
\begin{thebibliography}{10}

\bibitem{AS:handbook}
Milton Abramowitz and Irene~A. Stegun, \emph{Handbook of mathematical functions
  with formulas, graphs, and mathematical tables}, National Bureau of Standards
  Applied Mathematics Series, vol.~55, For sale by the Superintendent of
  Documents, U.S. Government Printing Office, Washington, D.C., 1964.

\bibitem{Archimedes:sphere}
{Archimedes of Syracuse}, \emph{On the sphere and cylinder}, ca. 225BC.

\bibitem{CE:new}
Henry Cohn and Noam Elkies, \emph{New upper bounds on sphere packings. {I}},
  Ann. of Math. (2) \textbf{157} (2003), no.~2, 689--714,
  \mbox{arXiv:math.MG/0110009}.

\bibitem{CS:splag}
John~H. Conway and Neil J.~A. Sloane, \emph{Sphere packings, lattices and
  groups}, 3rd ed., Grundlehren der mathematischen {Wissenschaften}, vol. 290,
  Springer-Verlag, New York, 1993.

\bibitem{CNR:smolyak}
R.~Cools, E.~Novak, and K.~Ritter, \emph{Smolyak's construction of cubature
  formulas of arbitrary trigonometric degree}, Computing \textbf{62} (1999),
  no.~2, 147--162.

\bibitem{CL:moderate}
Ronald Cools and James~N. Lyness, \emph{Three- and four-dimensional
  {$K$}-optimal lattice rules of moderate trigonometric degree}, Math. Comp.
  \textbf{70} (2001), no.~236, 1549--1567 (electronic).

\bibitem{CS:minimal}
Ronald Cools and Ian~H. Sloan, \emph{Minimal cubature formulae of trigonometric
  degree}, Math. Comp. \textbf{65} (1996), no.~216, 1583--1600.

\bibitem{Dehesa:transport}
Jes{\'u}s~S. Dehesa, \emph{Orthogonal polynomials in transport theories}, J.
  Phys. A \textbf{14} (1981), no.~2, 297--302.

\bibitem{Delsarte:linear}
P.~Delsarte, \emph{Bounds for unrestricted codes, by linear programming},
  Philips Res. Rep. \textbf{27} (1972), 272--289.

\bibitem{DGS:spherical}
P.~Delsarte, J.~M. Goethals, and J.~J. Seidel, \emph{Spherical codes and
  designs}, Geometriae Dedicata \textbf{6} (1977), no.~3, 363--388.

\bibitem{Frolov:lattice}
K.~K. Frolov, \emph{The connection of quadrature formulas and sublattices of
  the lattice of integer vectors}, Dokl. Akad. Nauk SSSR \textbf{232} (1977),
  no.~1, 40--43.

\bibitem{Ivanovic:formal}
I.~D. Ivanovi{\'c}, \emph{Formal state determination}, J. Math. Phys.
  \textbf{24} (1983), no.~5, 1199--1205.

\bibitem{KL:sphere}
G.~A. Kabatjanski\u{\i} and V.~I. Leven\v{s}te\u{\i}n, \emph{Bounds for
  packings on the sphere and in space}, Problemy Pereda\v ci Informacii
  \textbf{14} (1978), no.~1, 3--25.

\bibitem{Koike:spinor}
Kazuhiko Koike, \emph{Representations of spinor groups and the difference
  characters of {${\rm SO}(2n)$}}, Adv. Math. \textbf{128} (1997), no.~1,
  40--81.

\bibitem{Konig:cubature}
Hermann K{\"o}nig, \emph{Cubature formulas on spheres}, Advances in
  multivariate approximation (Witten-Bommerholz, 1998), Math. Res., vol. 107,
  Wiley-VCH, Berlin, 1999, pp.~201--211.

\bibitem{Krylov:calc}
Vladimir~Ivanovich Krylov, \emph{Approximate calculation of integrals}, The
  Macmillan Co., New York, 1962, Translated by Arthur H. Stroud.

\bibitem{Kuperberg:cubature}
Greg Kuperberg, \emph{Numerical cubature using error-correcting codes},
  \mbox{arXiv:math.NA/0402047}.

\bibitem{Moller:lower}
H.~Michael M\"oller, \emph{Lower bounds for the number of nodes in cubature
  formulae}, Numerische Integration (Tagung, Math. Forschungsinst.,
  Oberwolfach, 1978), Internat. Ser. Numer. Math., vol.~45, Birkh\"auser,
  Basel, 1979, pp.~221--230.

\bibitem{Mysovskikh:exact}
I.~P. Mysovskikh, \emph{Cubature formulas that are exact for trigonometric
  polynomials}, Dokl. Akad. Nauk SSSR \textbf{296} (1987), no.~1, 28--31.

\bibitem{NC:book}
Michael~A. Nielsen and Isaac~L. Chuang, \emph{Quantum computation and quantum
  information}, Cambridge University Press, Cambridge, 2000.

\bibitem{Noskov:periodic}
M.~V. Noskov, \emph{Formulas for the approximate integration of periodic
  functions}, Metody Vychisl. (1988), no.~15, 19--22, 178.

\bibitem{OS:unit}
A.~M. Odlyzko and N.~J.~A. Sloane, \emph{New bounds on the number of unit
  spheres that can touch a unit sphere in $n$ dimensions}, J. Combin. Theory
  Ser. A \textbf{26} (1979), no.~2, 210--214.

\bibitem{Rains:personal}
Eric Rains, 2003, personal communication.

\bibitem{Sidelnikov:7des}
V.~M. Sidelnikov, \emph{Spherical {$7$}-designs in {$2\sp n$}-dimensional
  {E}uclidean space}, J. Algebraic Combin. \textbf{10} (1999), no.~3, 279--288.

\bibitem{Sobolev:rotation}
S.~L. Sobolev, \emph{Cubature formulas on the sphere which are invariant under
  transformations of finite rotation groups}, Dokl. Akad. Nauk SSSR
  \textbf{146} (1962), 310--313.

\bibitem{Stanley:enumerative2}
Richard~P. Stanley, \emph{Enumerative combinatorics}, vol.~2, Cambridge
  University Press, 1999.

\bibitem{Stroud:some3}
A.~H. Stroud, \emph{Some approximate integration formulas of degree {$3$} for
  an {$n$}-dimensional simplex}, Numer. Math. \textbf{9} (1966), 38--45.

\bibitem{Stroud:calc}
\bysame, \emph{Approximate calculation of multiple integrals}, Prentice-Hall
  Inc., 1971, Prentice-Hall Series in Automatic Computation.

\bibitem{Stroud:more}
Arthur~H. Stroud, \emph{Quadrature methods for functions of more than one
  variable}, Ann. New York Acad. Sci. \textbf{86} (1960), 776--791 (1960).

\bibitem{Varadarajan:gtm}
V.~S. Varadarajan, \emph{Lie groups, {Lie} algebras, and their
  representations}, Graduate Texts in Mathematics, vol. 102, Springer-Verlag,
  New York, 1984, Reprint of the 1974 edition.

\bibitem{WX:triangle}
S.~Wandzura and H.~Xiao, \emph{Symmetric quadrature rules on a triangle},
  Comput. Math. Appl. \textbf{45} (2003), no.~12, 1829--1840.

\bibitem{WF:unbiased}
William~K. Wootters and Brian~D. Fields, \emph{Optimal state-determination by
  mutually unbiased measurements}, Ann. Physics \textbf{191} (1989), no.~2,
  363--381.

\bibitem{Xu:simplices}
Yuan Xu, \emph{Orthogonal polynomials and cubature formulae on spheres and on
  simplices}, Methods Appl. Anal. \textbf{5} (1998), no.~2, 169--184.

\bibitem{Xu:lower}
\bysame, \emph{Lower bound for the number of nodes of cubature formulae on the
  unit ball}, J. Complexity \textbf{19} (2003), no.~3, 392--402.

\bibitem{Magma}
\emph{Magma}, http://magma.maths.usyd.edu.au/.

\end{thebibliography}

\end{document}